# Natasha 2: Faster Non-Convex Optimization Than SGD
## — How to Swing By Saddle Points

(version 4)

Zeyuan Allen-Zhu
zeyuan@csail.mit.edu
Microsoft Research, Redmond

August 28, 2017*

**Abstract**

We design a stochastic algorithm to train any smooth neural network to $\varepsilon$-approximate local minima, using $O(\varepsilon^{-3.25})$ backpropagations. The best result was essentially $O(\varepsilon^{-4})$ by SGD. More broadly, it finds $\varepsilon$-approximate local minima of any smooth nonconvex function in rate $O(\varepsilon^{-3.25})$, with only oracle access to stochastic gradients.[1]

*V1 appeared on arXiv on this date. V2 and V3 polished writing. V4 was a deep revision and simplified proofs. This paper is built on, but should not be confused with, the *offline* method `Natasha1` [3] which only finds approximate stationary points.

[1]When this manuscript first appeared online, the best rate was $T = O(\varepsilon^{-4})$ by SGD. Several followups appeared after this paper. This includes SGD5 [5] and stochastic cubic Newton's method [46] both giving $T = O(\varepsilon^{-3.5})$, and `Neon+SCSG` [10, 48] which gives $T = O(\varepsilon^{-3.333})$. These rates are worse than $T = O(\varepsilon^{-3.25})$.

Our original method also requires oracle access to Hessian-vector products. However, this follow-up paper [10] enables us to replace the use of Hessian-vector products with stochastic gradient computations. We have revised this manuscript in V3 to reflect this change.

# 1 Introduction

In diverse world of deep learning research has given rise to numerous architectures for neural networks (convolutional ones, long short term memory ones, etc). However, to this date, the underlying training algorithms for neural networks are still stochastic gradient descent (SGD) and its heuristic variants. In this paper, we address the problem of designing a new algorithm that has provably faster running time than the best known result for SGD.

Mathematically, we study the problem of online stochastic nonconvex optimization:

$$\min_{x\in\mathbb{R}^d} \Big\{ f(x) \stackrel{\text{def}}{=} \mathbb{E}_i[f_i(x)] = \tfrac{1}{n}\textstyle\sum_{i=1}^{n} f_i(x) \Big\} \tag{1.1}$$

where both $f(\cdot)$ and each $f_i(\cdot)$ can be nonconvex. We want to study

*online algorithms* to find approximate *local minimum* of $f(x)$.

Here, we say an algorithm is online if its complexity is independent of $n$. This tackles the big-data scenarios when $n$ is extremely large or even infinite.[2]

Nonconvex optimization arises prominently in large-scale machine learning. Most notably, training deep neural networks corresponds to minimizing $f(x)$ of this average structure: each training sample $i$ corresponds to one loss function $f_i(\cdot)$ in the summation. This average structure allows one to perform stochastic gradient descent (SGD) which uses a random $\nabla f_i(x)$ —corresponding to computing backpropagation once— to approximate $\nabla f(x)$ and performs descent updates.

The standard goal of nonconvex optimization with provable guarantee is to find *approximate local minima*. This is not only because finding the *global* one is NP-hard, but also because there exist rich literature on *heuristics* for turning a local-minima finding algorithm into a global one. This includes random seeding, graduated optimization [25] and others. Therefore, faster algorithms for finding approximate local minima translate into faster *heuristic* algorithms for finding global minimum.

On a separate note, experiments [16, 17, 24] suggest that fast convergence to approximate local minima may be sufficient for training neural nets, while convergence to stationary points (i.e., points that may be saddle points) is *not*. In other words, we need to *avoid saddle points.*

## 1.1 Classical Approach: Escaping Saddle Points Using Random Perturbation

One natural way to avoid saddle points is to use randomness to escape from it, whenever we meet one. For instance, Ge et al. [22] showed, by injecting random perturbation, SGD will not be stuck in saddle points: whenever SGD moves into a saddle point, randomness shall help it escape. This partially explains why SGD performs well in deep learning.[3] Jin et al. [27] showed, equipped with random perturbation, full gradient descent (GD) also escapes from saddle points. Being easy to implement, however, we raise two main efficiency issues regarding this classical approach:

- Issue 1. If we want to escape from saddle points, is random perturbation the only way? Moving in a random direction is "blind" to the Hessian information of the function, and thus can we escape from saddle points faster?
- Issue 2. If we want to avoid saddle points, is it really necessary to first move close to saddle points and then *escape* from them? Can we design an algorithm that can somehow avoid saddle points without ever moving close to them?

[2]All of our results in this paper apply to the case when $n$ is infinite, meaning $f(x) = \mathbb{E}_i[f_i(x)]$, because we focus on *online* methods. However, we still introduce $n$ to simplify notations.

[3]In practice, stochastic gradients naturally incur "random noise" and adding perturbation may not be needed.

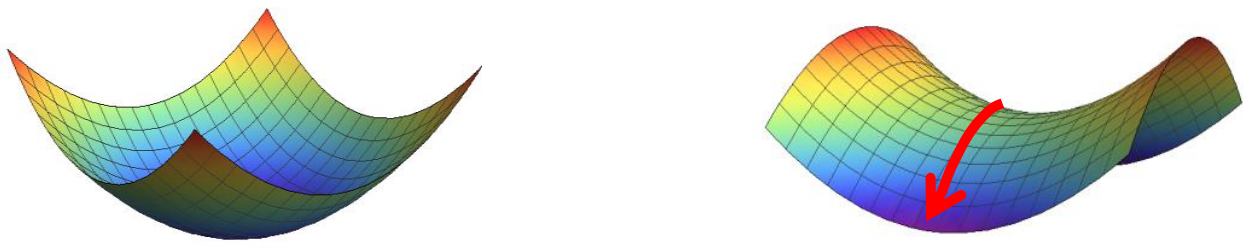

Figure 1: Local minimum (left), saddle point (right) and its negative-curvature direction.

### 1.2 Our Resolutions

**Resolution to Issue 1: Efficient Use of Hessian.** Mathematically, instead of using a random perturbation, the negative eigenvector of $\nabla^2 f(x)$ (a.k.a. the negative-curvature direction of $f(\cdot)$ at $x$) gives us a *better* direction to escape from saddle points. See Figure 1.

To make it concrete, suppose we apply power method on $\nabla^2 f(x)$ to find its most negative eigenvector. [4]If we run power method for 0 iteration, then it gives us a totally random direction; if we run it for more iterations, then it converges to the most negative eigenvector of $\nabla^2 f(x)$. Unfortunately, applying power method is unrealistic because $\nabla^2 f(x)$ is a large matrix and $f(x) = \frac{1}{n}\sum_i f_i(x)$ can possibly have infinite pieces.

We propose to use Oja's algorithm [37] to approximate power method. Oja's algorithm can be viewed as an online variant of power method, and requires only (stochastic) matrix-vector product computations. In our setting, this is the same as (stochastic) Hessian-vector products —namely, computing $\nabla^2 f_i(x) \cdot w$ for arbitrary vectors $w \in \mathbb{R}^d$ and random indices $i \in [n]$. It is a known fact that computing Hessian-vector products is as cheap as computing stochastic gradients (see Remark 1.1), and thus we can use Oja's algorithm to escape from saddle points. (This requires the recent convergence analysis of Oja's algorithm by Allen-Zhu and Li [9].)

*Remark* 1.1. Computing (stochastic) Hessian-vector products is as cheap as computing stochastic gradients. This can be seen in at least two ways.

- If $f_i(x)$ is described by a size-$S$ arithmetic circuit, then computing $\nabla f_i(x)$ and $\nabla^2 f_i(x) \cdot w$ both cost running time $O(S)$ due to the chain rule of derivative [38]. For *training neural networks*, computing $\nabla f_i(x)$ requires one backpropagation; but $\nabla^2 f_i(x) \cdot w$ can also be implemented via one backpropagation, for a network of roughly the same size [38, 42]. In practice, some reported that $\nabla^2 f_i(x) \cdot w$ is twice expensive to compute as $\nabla f_i(x)$ [42] in training neural networks.
- One can also use $\frac{\nabla f_i(x+qw) - \nabla f_i(x)}{q}$ to approximate $\nabla^2 f_i(x) \cdot w$ when $q$ is a small positive constant. This is not only used in practice, but can also be made mathematically rigorous for certain algorithms, including the one we shall introduce in this paper.[5]

**Resolution to Issue 2: Swing by Saddle Points.** If the function is sufficiently smooth,[6] then any point close to a saddle point must have a negative curvature. Therefore, as long as we are close to saddle points, we can already use Oja's algorithm to find such negative curvature, and move in its direction to decrease the objective, see Figure 2(a).

Therefore, we are left only with the case that point is not close to any saddle point. Using smoothness of $f(\cdot)$, this gives a "safe zone" near the current point, in which there is no strict saddle point, see Figure 2(b). Intuitively, we wish to use the property of safe zone to design an algorithm that decreases the objective faster than SGD. Formally, $f(\cdot)$ inside this safe zone must be of "bounded nonconvexity," meaning that its eigenvalues of the Hessian are always greater than some negative threshold $-\sigma$ (where $\sigma$ depends on how long we run Oja's algorithm). Intuitively, the greater $\sigma$ is, then the more non-convex $f(x)$ is. We wish to design an (online) stochastic first-order method whose running time scales with $\sigma$.

[4]In our context, up to normalization, power method outputs a vector $v' = (\mathbf{I} - \eta\nabla^2 f(x))^M \cdot v$ where $v$ is a random vector, $\eta > 0$ is some small learning rate, and $M \geq 0$ is the number of iterations.

[5]In a follow-up work, Allen-Zhu and Li [10] showed that Hessian-vector product computations in this paper can be replaced by $\frac{\nabla f_i(x+qw) - \nabla f_i(x)}{q}$. We discuss this more in Section 5.1.

[6]As we shall see, smoothness is necessary for finding approximate local minima with provable guarantees.

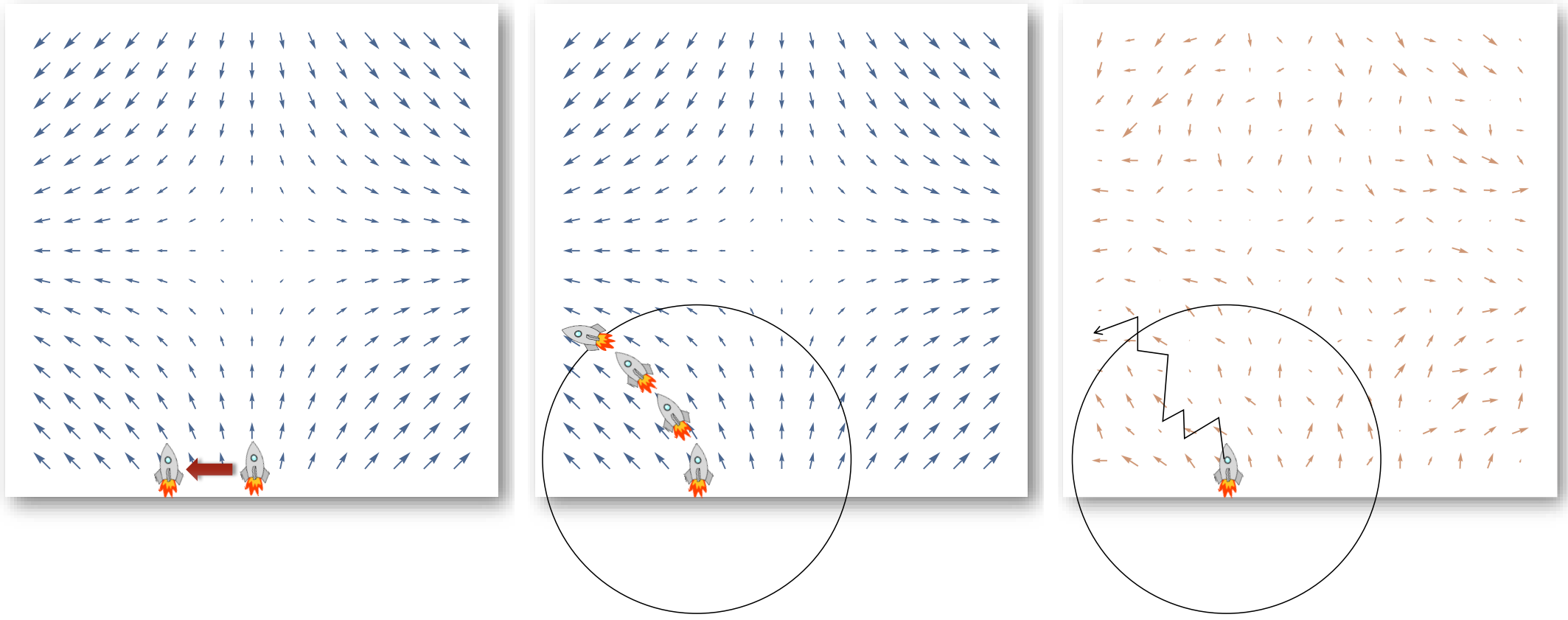

(a) (b) (c)

Figure 2: Illustration of $\texttt{Natasha2}^{\mathsf{full}}$ — how to swing by a saddle point.
(a) move in a negative curvature direction if there is any (by applying Oja's algorithm)
(b) swing by a saddle point without entering its neighborhood (wishful thinking)
(c) swing by a saddle point using only stochastic gradients (by applying $\texttt{Natasha1.5}^{\mathsf{full}}$)

Unfortunately, classical stochastic methods such as SGD or SCSG [30] cannot make use of this nonconvexity parameter $\sigma$. The only known ones that can make use of $\sigma$ are offline algorithms. In this paper, we design a new stochastic first-order method `Natasha1.5`

**Theorem 1** (informal)**.** `Natasha1.5` *finds* $x$ *with* $\|\nabla f(x)\| \leq \varepsilon$ *in rate* $T = O\big(\frac{1}{\varepsilon^3} + \frac{\sigma^{1/3}}{\varepsilon^{10/3}}\big)$ .

Finally, we put `Natasha1.5` together with Oja's to construct our final algorithm `Natasha2`:

**Theorem 2** (informal)**.** `Natasha2` *finds* $x$ *with* $\|\nabla f(x)\| \leq \varepsilon$ *and* $\nabla^2 f(x) \succeq -\delta \mathbf{I}$ *in rate* $T = \widetilde{O}\big(\frac{1}{\delta^5} + \frac{1}{\delta\varepsilon^3} + \frac{1}{\varepsilon^{3.25}}\big)$*. In particular, when* $\delta \geq \varepsilon^{1/4}$*, this gives* $T = \widetilde{O}\big(\frac{1}{\varepsilon^{3.25}}\big)$*.*

In contrast, the convergence rate of SGD was $T = \widetilde{O}(\mathsf{poly}(d) \cdot \varepsilon^{-4})$ [22].

### 1.3 Follow-Up Results

Since the original appearance of this work, there has been a lot of progress in stochastic nonconvex optimization. Most notably,

- If one *swings by* saddle points using Oja's algorithm and SGD variants (instead of `Natasha1.5`), the convergence rate is $T = \widetilde{O}(\varepsilon^{-3.5})$ [5].
- If one applies SGD and only *escapes from* saddle points using Oja's algorithm, the convergence rate is $T = \widetilde{O}(\varepsilon^{-4})$ [10, 48].
- If one applies SCSG and only *escapes from* saddle points using Oja's algorithm, the convergence rate is $T = \widetilde{O}(\varepsilon^{-3.333})$ [10, 48].
- If one applies a stochastic version of cubic regularized Newton's method, the convergence rate is $T = \widetilde{O}(\varepsilon^{-3.5})$ [46].
- If $f(x)$ is of $\sigma$-bounded nonconvexity, the SGD4 method [5] gives rate $T = \widetilde{O}(\varepsilon^{-2} + \sigma\varepsilon^{-4})$.

We include these results in Table 1 for a close comparison.

| | algorithm | gradient complexity $T$ | variance bound | Lipschitz smooth | 2nd-order smooth |
|---|---|---|---|---|---|
| convex only | SGD1 [5, 23] | $O(\varepsilon^{-2.667})$ | needed | needed | no |
| | SGD2 [5] | $O(\varepsilon^{-2.5})$ ♯ | needed | needed | no |
| | SGD3 [5] | $\widetilde{O}(\varepsilon^{-2})$ ♯ | needed | needed | no |
| approximate stationary points | SGD (folklore) | $O(\varepsilon^{-4})$ (see Appendix B) | needed | needed | no |
| | SCSG [30] | $O(\varepsilon^{-3.333})$ | needed | needed | no |
| | Natasha1.5 | $O(\varepsilon^{-3} + \sigma^{1/3}\varepsilon^{-3.333})$ (see Theorem 1) | needed | needed | no |
| | SGD4 [5] | $\widetilde{O}(\varepsilon^{-2} + \sigma\varepsilon^{-4})$ ♯ | needed | needed | no |
| approximate local minima | perturbed SGD [22] | $\widetilde{O}(\varepsilon^{-4} \cdot \mathsf{poly}(d))$ | needed | needed | needed |
| | Natasha2 | $\widetilde{O}(\varepsilon^{-3.25})$ (see Theorem 2) | needed | needed | needed |
| | NEON + SGD [10, 48] | $\widetilde{O}(\varepsilon^{-4})$ ♯ | needed | needed | needed |
| | cubic Newton [46] | $\widetilde{O}(\varepsilon^{-3.5})$ ♯ | needed | needed | needed |
| | SGD5 [5] | $\widetilde{O}(\varepsilon^{-3.5})$ ♯ | needed | needed | needed |
| | NEON + SCSG [10, 48] | $\widetilde{O}(\varepsilon^{-3.333})$ ♯ | needed | needed | needed |

Table 1: Comparison of ***online*** methods for finding $\|\nabla f(x)\| \le \varepsilon$. Following tradition, in these complexity bounds, we assume variance and smoothness parameters as constants, and only show the dependency on $n, d, \varepsilon$ and the bounded nonconvexity parameter $\sigma \in (0, 1)$. We use ♯ to indicate results that appeared after this paper.

**Remark 1.** Variance bounds must be needed for online methods.

**Remark 2.** Lipschitz smoothness must be needed for achieving even approximate stationary points.

**Remark 3.** Second-order smoothness must be needed for achieving approximate local minima.

### 1.4 Roadmap

We introduce necessary notations in Section 2, and give high-level intuitions and pseudocodes of Natasha1.5 and Natasha2 respectively in Section 3 and Section 4. In Section 5, we review Oja's algorithm and prove some auxiliary theorems. In Section 6 and Section 7 we give full proofs to Natasha1.5 and Natasha2. Some more related work is discussed in Section A, and proofs for SGD and GD for finding approximate stationary points are included in Section B for completeness' sake.

## 2 Preliminaries

Throughout this paper, we denote by $\|\cdot\|$ the Euclidean norm. We use $i \in_R [n]$ to denote that $i$ is generated from $[n] = \{1, 2, \dots, n\}$ uniformly at random. We denote by $\nabla f(x)$ the gradient of function $f$ if it is differentiable, and $\partial f(x)$ any subgradient if $f$ is only Lipschitz continuous. We denote by $\mathbb{I}[event]$ the indicator function of probabilistic events.

We denote by $\|\mathbf{A}\|_2$ the spectral norm of matrix $\mathbf{A}$. For symmetric matrices $\mathbf{A}$ and $\mathbf{B}$, we write $\mathbf{A} \succeq \mathbf{B}$ to indicate that $\mathbf{A} - \mathbf{B}$ is positive semidefinite (PSD). Therefore, $\mathbf{A} \succeq -\sigma\mathbf{I}$ if and only if all eigenvalues of $\mathbf{A}$ are no less than $-\sigma$. We denote by $\lambda_{\min}(\mathbf{A})$ and $\lambda_{\max}(\mathbf{A})$ the minimum and maximum eigenvalue of a symmetric matrix $\mathbf{A}$.

Recall some definitions on strong convexity (SC), bounded nonconvexity, and smoothness.

**Definition 2.1.** *For a function* $f\colon \mathbb{R}^d \to \mathbb{R}$,

- $f$ *is* $\sigma$*-strongly convex if* $\forall x, y \in \mathbb{R}^d$, *it satisfies* $f(y) \ge f(x) + \langle \partial f(x), y - x\rangle + \frac{\sigma}{2}\|x - y\|^2$.

- *$f$ is of $\sigma$-bounded nonconvexity (or $\sigma$-**nonconvex** for short) if $\forall x, y \in \mathbb{R}^d$, it satisfies $f(y) \geq f(x) + \langle \partial f(x), y - x\rangle - \frac{\sigma}{2}\|x-y\|^2$.* [7]
- *$f$ is $L$-Lipschitz smooth (or $L$-**smooth** for short) if $\forall x, y \in \mathbb{R}^d$, $\|\nabla f(x) - \nabla f(y)\| \leq L\|x-y\|$.*
- *$f$ is second-order $L_2$-Lipschitz smooth (or $L_2$-**second-order smooth** for short) if $\forall x, y \in \mathbb{R}^d$, it satisfies $\|\nabla^2 f(x) - \nabla^2 f(y)\|_2 \leq L_2\|x-y\|$.*

These definitions have other equivalent forms, see textbook [33].

**Definition 2.2.** *For composite function $F(x) = \psi(x) + f(x)$ where $\psi(x)$ is proper convex, given a parameter $\eta > 0$, the **gradient mapping** of $F(\cdot)$ at point $x$ is*

$$\mathcal{G}_{F,\eta}(x) \stackrel{\text{def}}{=} \frac{1}{\eta}\big(x - x'\big) \qquad \text{where} \qquad x' = \arg\min_y \big\{\psi(y) + \langle \nabla f(x), y\rangle + \frac{1}{2\eta}\|y-x\|^2\big\}$$

*In particular, if $\psi(\cdot) \equiv 0$, then $\mathcal{G}_{F,\eta}(x) \equiv \nabla f(x)$.*

## 3 Natasha 1.5: Finding Approximate Stationary Points

We first make a detour to study how to find approximate stationary points using only first-order information. A point $x \in \mathbb{R}^d$ is an $\varepsilon$-approximate stationary point[8] of $f(x)$ if it satisfies $\|\nabla f(x)\| \leq \varepsilon$. Let *gradient complexity* $T$ be the number of computations of $\nabla f_i(x)$.

Before 2015, nonconvex first-order methods give rise to two convergence rates. SGD converges in $T = O(\varepsilon^{-4})$ and GD converges $T = O(n\varepsilon^{-2})$. The proofs of both are simple (see Appendix B for completeness). In particular, the convergence of SGD relies on two *minimal* assumptions

$$f(x) \text{ has bounded variance } \mathcal{V}, \text{ meaning } \mathbb{E}_i[\|\nabla f_i(x) - \nabla f(x)\|^2] \leq \mathcal{V}, \text{ and} \tag{A1}$$

$$f(x) \text{ is } L\text{-Lipschitz smooth, meaning } \|\nabla f(x) - \nabla f(y)\| \leq L \cdot \|x-y\|. \tag{A2'}$$

*Remark* 3.1. Both assumptions are necessary to design online algorithms for finding stationary points.[9] For *offline* algorithms —like GD— the first assumption is not needed.

Since 2016, the convergence rates have been improved to $T = O(n + n^{2/3}\varepsilon^{-2})$ for offline methods [6, 39], and to $T = O(\varepsilon^{-10/3})$ for online algorithms [30]. Both results are based on the SVRG (stochastic variance reduced gradient) method, and assume additionally (note (A2) implies (A2'))

$$\text{each } f_i(x) \text{ is } L\text{-Lipschitz smooth.} \tag{A2}$$

Lei et al. [30] gave their algorithm a new name, SCSG (stochastically controlled stochastic gradient).

**Bounded Non-Convexity.** In recent works [3, 13], it has been proposed to study a more refined convergence rate, by assuming that $f(x)$ is of *$\sigma$-bounded nonconvexity* (or $\sigma$-nonconvex), meaning

$$\text{all the eigenvalues of } \nabla^2 f(x) \text{ lie in } [-\sigma, L] \tag{A3}$$

for some $\sigma \in (0, L]$. This parameter $\sigma$ is analogous to the strong-convexity parameter $\mu$ in convex optimization, where all the eigenvalues of $\nabla^2 f(x)$ lie in $[\mu, L]$ for some $\mu > 0$.

[7] Previous authors also refer to this notion as "approximate convex", "almost convex", "hypo-convex", "semi-convex", or "weakly-convex." We call it $\sigma$-nonconvex to stress the point that $\sigma$ can be as large as $L$ (any $L$-smooth function is automatically $L$-nonconvex).

[8] Historically, in first-order literatures, $x$ is called $\varepsilon$-approximate if $\|\nabla f(x)\|^2 \leq \varepsilon$; in second-order literatures, $x$ is $\varepsilon$-approximate if $\|\nabla f(x)\| \leq \varepsilon$. We adapt the latter notion following Polyak and Nesterov [34, 36].

[9] For instance, if the variance $\mathcal{V}$ is unbounded, we cannot even tell if a point $x$ satisfies $\|\nabla f(x)\| \leq \varepsilon$ using finite samples. Also, if $f(x)$ is not Lipschitz smooth, it may contain sharp turning points (e.g., behaves like absolute value function $|x|$); in this case, finding $\|\nabla f(x)\| \leq \varepsilon$ can be as hard as finding $\|\nabla f(x)\| = 0$, and is NP-hard in general.

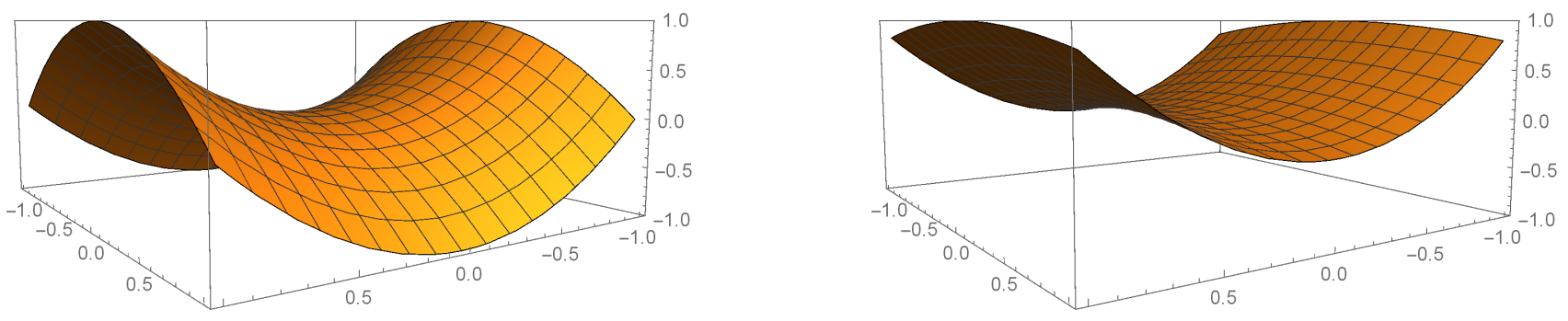

(a) 1-smooth & 1-nonconvex (b) 1-smooth & 0.2-nonconvex

Figure 3: Nonconvex functions with bounded nonconvexity.

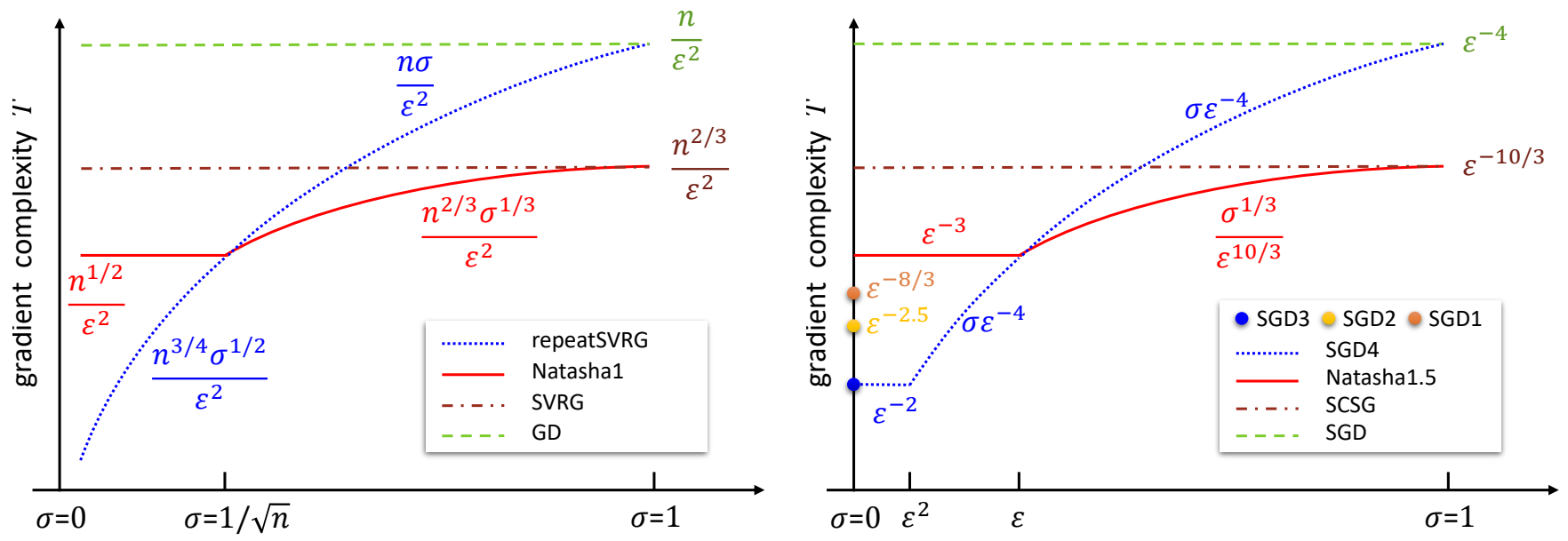

(a) offline first-order methods (b) *online* first-order methods

Figure 4: Comparison of first-order methods for finding $\varepsilon$-approximate stationary points of a $\sigma$-nonconvex function. For simplicity, in the plots we let $L = 1$ and $\mathcal{V} = 1$. The results SGD2/3/4 appeared after this work.

As examples, Figure 3(a) is 1-nonconvex and Figure 3(b) is 0.2-nonconvex. In our illustrative process to "swing by a saddle point," the function inside safe zone —see Figure 2(b)— is also of bounded nonconvexity. Since larger $\sigma$ means the function is "more non-convex" and thus harder to optimize, can we design algorithms with gradient complexity $T$ as an *increasing function of* $\sigma$ ?

*Remark* 3.2. Most methods (SGD, SCSG, SVRG and GD) do not run faster if $\sigma < L$, at least in theory. More work is thus needed.

In the *offline* setting, two methods are known to make use of parameter $\sigma$. One is `repeatSVRG`, implicitly in [13] and formally in [3]. The other is `Natasha1` [3]. `repeatSVRG` performs better when $\sigma \le L/\sqrt{n}$ and `Natasha1` performs better when $\sigma \ge L/\sqrt{n}$. See Figure 4(a) and Table 2.

Before this work, no online method is known to take advantage of $\sigma$.

## 3.1 Our Theorem

We show that, under (A1), (A2) and (A3), one can non-trivially extend `Natasha1` to an *online* version, taking advantage of $\sigma$, and achieving better complexity than `SCSG`.

Let $\Delta_f$ be any upper bound on $f(x_0) - f(x^*)$ where $x_0$ is the starting point. In this section, to present the simplest results, we use the big-$O$ notion to hide dependency in $\Delta_f$ and $\mathcal{V}$. In Section 6, we shall add back such dependency and as well as support the existence of a proximal term. (That is, to minimize $\psi(x) + f(x)$ where $\psi(x)$ is a proper convex simple function.)

Under such simplified notations, our main theorem can be stated as follows.

**Theorem 1** (simple). *Under (A1), (A2) and (A3), using the big-O notion to hide dependency in* $\Delta_f$ *and* $\mathcal{V}$*, we have for every* $\varepsilon \in (0, \frac{\sigma}{L}]$*, letting*

$$B = \Theta\left(\frac{1}{\varepsilon^2}\right) \quad , \quad T = \Theta\left(\frac{L^{2/3}\sigma^{1/3}}{\varepsilon^{10/3}}\right) \quad \textit{and} \quad \alpha = \Theta\left(\frac{\varepsilon^{4/3}}{\sigma^{1/3}L^{2/3}}\right)$$

**Algorithm 1** $\texttt{Natasha1.5}(F, x^{\varnothing}, B, T', \alpha)$

**Input:** $f(\cdot) = \frac{1}{n}\sum_{i=1}^{n} f_i(x)$, starting vector $x^{\varnothing}$, epoch length $B \in [n]$, epoch count $T' \geq 1$, learning rate $\alpha > 0$.

$\widehat{\mathsf{x}} \leftarrow x^{\varnothing}$; $\quad p \leftarrow \Theta((\sigma/\varepsilon L)^{2/3})$; $\quad m \leftarrow B/p$; $\quad X \leftarrow []$;
**for** $k \leftarrow 1$ **to** $T'$ **do** ◇ *$T'$ epochs each of length $B$*
  $\widetilde{\mathbf{x}} \leftarrow \widehat{\mathsf{x}}$; $\mu \leftarrow \frac{1}{B}\sum_{i\in S}\nabla f_i(\widetilde{\mathbf{x}})$ where $S$ is a uniform random subset of $[n]$ with $|S| = B$;
  **for** $s \leftarrow 0$ **to** $p-1$ **do** ◇ *$p$ sub-epochs each of length $m$*
    $x_0 \leftarrow \widehat{\mathsf{x}}$; $\quad X \leftarrow [X, \widehat{\mathsf{x}}]$;
    **for** $t \leftarrow 0$ **to** $m-1$ **do**
      $\widetilde{\nabla} \leftarrow \nabla f_i(x_t) - \nabla f_i(\widetilde{\mathbf{x}}) + \mu + 2\sigma(x_t - \widehat{\mathsf{x}})$ where $i \in_R [n]$
      $x_{t+1} = x_t - \alpha\widetilde{\nabla}$;
    **end for**
    $\widehat{\mathsf{x}} \leftarrow$ a random choice from $\{x_0, x_1, \dots, x_{m-1}\}$; ◇ *in practice, choose the average*
  **end for**
**end for**
$\widehat{\mathsf{y}} \leftarrow$ a random vector in $X$. ◇ *in practice, simply return $\widehat{\mathsf{y}}$*
$g(x) \stackrel{\text{def}}{=} f(x) + \sigma\|x - \widehat{\mathsf{y}}\|^2$ and use convex SGD to minimize $g(x)$ for $T_{\mathsf{sgd}} = T'B$ iterations.
**return** $x^{\mathsf{out}} \leftarrow$ the output of SGD.

*we have that* $\texttt{Natasha1.5}(f, x^{\varnothing}, B, T/B, \alpha)$ *outputs a point* $x^{\mathsf{out}}$ *with* $\mathbb{E}[\|\nabla f(x^{\mathsf{out}})\|] \leq \varepsilon$, *and needs* $O(T)$ *computations of stochastic gradients. (See also Figure 4(b).)*

We emphasize that the additional factor $\sigma^{1/3}$ in the numerator of $T$ shall become our key to achieve faster algorithm for finding approximate local minima in Section 4. Also, if the requirement $\varepsilon \leq \frac{\sigma}{L}$ is not satisfied, one can replace $\sigma$ with $\varepsilon L$; accordingly, $T$ becomes $O\big(\frac{L}{\varepsilon^3} + \frac{L^{2/3}\sigma^{1/3}}{\varepsilon^{10/3}}\big)$

We note that the SGD4 method of [5] (which appeared after this paper) achieves $T = O\big(\frac{L}{\varepsilon^2} + \frac{\sigma}{\varepsilon^4}\big)$. It is better than Natasha1.5 only when $\sigma \leq \varepsilon L$. We compare them in Figure 4(b), and emphasize that it is necessary to use Natasha1.5 (rather than SGD4) to design Natasha2 of the next section.

**Extension.** In fact, we show Theorem 1 in a more general *proximal* setting. That is, to minimize $F(x) \stackrel{\text{def}}{=} f(x) + \psi(x)$ where $\psi(x)$ is proper convex function that can be *non-smooth*. For instance, if $\psi(x)$ is the indicator function of a convex set, then Problem (1.1) becomes constraint minimization; and if $\psi(x) = \|x\|_1$, we encourage sparsity. At a first reading of its proof, one can assume $\psi(x) \equiv 0$.

## 3.2 Our Intuition

We first recall the main idea of the SVRG method [28, 50], which is an *offline* algorithm. SVRG divides iterations into epochs, each of length $n$. It maintains a snapshot point $\widetilde{\mathbf{x}}$ for each epoch, and computes the full gradient $\nabla f(\widetilde{\mathbf{x}})$ only for snapshots. Then, in each iteration $t$ at point $x_t$, SVRG defines gradient estimator $\widetilde{\nabla} f(x_t) \stackrel{\text{def}}{=} \nabla f_i(x_t) - \nabla f_i(\widetilde{\mathbf{x}}) + \nabla f(\widetilde{\mathbf{x}})$ which satisfies $\mathbb{E}_i[\widetilde{\nabla} f(x_t)] = \nabla f(x_t)$, and performs proximal update $x_{t+1} \leftarrow x_t - \alpha\widetilde{\nabla} f(x_t)$ for learning rate $\alpha$.

For minimizing non-convex functions, SVRG does not take advantage of parameter $\sigma$ *even if* the learning rate can be adapted to $\sigma$. This is because SVRG (and in fact SGD and GD too) rely on gradient-descent analysis to argue for objective decrease *per iteration*. This is blind to $\sigma$.[10]

[10]These results argue for objective decrease per iteration, of the form $f(x_t) - f(x_{t+1}) \geq \frac{\alpha}{2}\|\nabla f(x_t)\|^2 - \frac{\alpha^2 L}{2}\mathbb{E}\big[\|\nabla f(x_t) - \widetilde{\nabla} f(x_t)\|^2\big]$. Unlike mirror-descent analysis, this inequality cannot take advantage of the bounded nonconvexity parameter of $f(x)$. For readers interested in the difference between gradient and mirror descent, see [11].

The prior work Natasha1 takes advantage of $\sigma$. Natasha1 is similar to SVRG, but it further divides each epoch into sub-epochs, each with a starting vector $\widehat{\mathsf{x}}$. Then, it replaces $\widetilde{\nabla} f(x_t)$ with $\widetilde{\nabla} f(x_t) + 2\sigma(x_t - \widehat{\mathsf{x}})$. This is equivalent to replacing $f(x)$ with $f(x) + \sigma\|x - \widehat{\mathsf{x}}\|^2$, where the center $\widehat{\mathsf{x}}$ changes every sub-epoch. We view this additional term $2\sigma(x_t - \widehat{\mathsf{x}})$ as a type of *retraction*. Conceptually, it stabilizes the algorithm by moving a bit in the backward direction. Technically, it enables us to perform only mirror-descent type of analysis, and thus bypass the issue of SVRG.

**Our Algorithm.** Both SVRG and Natasha1 are offline methods, because the gradient estimator requires the full gradient computation $\nabla f(\widetilde{\mathbf{x}})$ at snapshots $\widetilde{\mathbf{x}}$. A natural fix —originally studied by practitioners but first formally analyzed by Lei et al. [30]— is to replace the computation of $\nabla f(\widetilde{\mathbf{x}})$ with $\frac{1}{|S|}\sum_{i\in S} \nabla f_i(\widetilde{\mathbf{x}})$, for a random batch $S \subseteq [n]$ with fixed cardinality $B := |S| \ll n$. This allows us to shorten the epoch length from $n$ to $B$, thus turning SVRG and Natasha1 into *online* methods.

How large should we pick $B$? By Chernoff bound, we wish $B \approx \frac{1}{\varepsilon^2}$ because our desired accuracy is $\varepsilon$. One can thus *hope* to replace the parameter $n$ in the complexities of SVRG and Natasha1.5 (ignoring the dependency on $L$):

$$T = O\big(n + n^{2/3}\varepsilon^{-2}\big) \quad \text{and} \quad T = O\big(n + n^{1/2}\varepsilon^{-2} + \sigma^{1/3} n^{2/3}\varepsilon^{-2}\big)$$

with $B \approx \frac{1}{\varepsilon^2}$. This "wishful thinking" gives

$$T = O\big(\varepsilon^{-10/3}\big) \quad \text{and} \quad T = O\big(\varepsilon^{-3} + \sigma^{1/3}\varepsilon^{-10/3}\big).$$

These are exactly the results achieved by SCSG [30] and to be achieved by our new Natasha1.5.

*Unfortunately,* Chernoff bound itself is not sufficient in getting such rates. Let

$$\mathbf{e} \stackrel{\text{def}}{=} \tfrac{1}{|S|}\textstyle\sum_{i\in S} \nabla f_i(\widetilde{\mathbf{x}}) - \nabla f(\widetilde{\mathbf{x}})$$

denote the bias of this new gradient estimator, then when performing iterative updates, this bias $\mathbf{e}$ gives rise to two types of error terms: "first-order error" terms —of the form $\langle \mathbf{e}, x - y\rangle$— and "second-order error" term $\|\mathbf{e}\|^2$. Chernoff bound ensures that the second-order error $\mathbb{E}_S[\|\mathbf{e}\|^2] \leq \varepsilon^2$ is bounded. However, first-order error terms are the true bottlenecks.

In the offline method SCSG, Lei et al. [30] carefully performed updates so that all "first-order errors" cancel out. To the best of our knowledge, this analysis cannot take advantage of $\sigma$ even if the algorithm knows $\sigma$. (Again, for experts, this is because SCSG is based on gradient-descent type of analysis but not mirror-descent.)

In Natasha1.5, we use the aforementioned retraction to ensure that all points in a single sub-epoch are close to each other (based on mirror-descent type of analysis). Then, we use Young's inequality to bound $\langle \mathbf{e}, x - y\rangle$ by $\frac{1}{2}\|\mathbf{e}\|^2 + \frac{1}{2}\|x - y\|^2$. In this equation, $\|\mathbf{e}\|^2$ is already bounded by Chernoff concentration, and $\|x - y\|^2$ can also be bounded as long as $x$ and $y$ are within the same sub-epoch. This summarizes the high-level technical contribution of Natasha1.5.

We formally state Natasha1.5 in Algorithm 1, and it uses big-$O$ notions to hide dependency in $L$, $\Delta_f$, and $\mathcal{V}$. The more general code to take care of the proximal term is in Algorithm 3 of Section 6.

*Remark* 3.3. The SCSG method by Lei et al. [30] is in fact SVRG plus two modifications. The first is to reduce $n$ to $B$ as discussed above. The second is to *randomly* stop an epoch so that its length forms a memoryless geometric distribution. They call this algorithm SCSG. As we have demonstrated in this paper, this random stopping technique is not really necessary.

*Remark* 3.4. In our pseudocode of Natasha1.5, we have twice selected *random points* in the computation history. We use a random point within each subepoch $\{x_0, x_1, \ldots, x_m\}$ as the next starting point $\widehat{\mathsf{x}}$, and select $\widehat{\mathsf{y}}$ as a random copy of $\widehat{\mathsf{x}}$. Selecting random points is necessary for theoretical purpose (and was even present in the basic proofs of non-convex GD and SGD, see Appendix B), but we recommend selecting the *last points* in practice.

*Remark* 3.5. In our pseudocode of `Natasha1.5`, we have an addition pruning step in Line 14. That is, instead of directly outputting $\widehat{\mathsf{y}}$, it regularizes the function $f(x) + \sigma\|x - \widehat{\mathsf{y}}\|^2$ to make it convex, and then apply SGD to minimize it to some sufficient accuracy. This pruning step is also for the purpose of proving theoretical convergence, and is not necessary in practice.

# 4 Natasha 2: Finding Approximate Local Minima

Stochastic gradient descent (SGD) find approximate local minima [22], under (A1), (A2) and an additional assumption (A4):

$$f(x) \text{ is second-order } L_2\text{-Lipschitz smooth, meaning } \|\nabla^2 f(x) - \nabla^2 f(y)\|_2 \le L_2 \cdot \|x - y\|. \quad \text{(A4)}$$

*Remark* 4.1. (A4) is necessary to make the task of find approximate local minima meaningful, for the same reason Lipschitz smoothness was needed for finding stationary points.

**Definition 4.2.** *We say $x$ is an $(\varepsilon, \delta)$-approximate local minimum of $f(x)$ if*[11]

$$\|\nabla f(x)\| \le \varepsilon \quad \textit{and} \quad \nabla^2 f(x) \succeq -\delta \mathbf{I} \ ,$$

*or $\varepsilon$-approximate local minimum if it is $(\varepsilon, \varepsilon^{1/C})$-approximate local minimum for constant $C \ge 1$.*

(Approximate local minima are also known as approximate second-order critical points, but may not be close to *any* exact local minimum.)

Before our work, Ge et al. [22] is the only result that gives provable *online* complexity for finding approximate local minima. Other previous results, including `SVRG`, `SCSG`, `Natasha1`, and even `Natasha1.5`, do not find approximate local minima and may be stuck at saddle points.[12] Ge et al. [22] showed that, hiding factors that depend on $L$, $L_2$ and $\mathcal{V}$, SGD finds an $\varepsilon$-approximate local minimum of $f(x)$ in gradient complexity $T = O(\mathsf{poly}(d)\varepsilon^{-4})$. This $\varepsilon^{-4}$ factor seems necessary since SGD needs $T \ge \Omega(\varepsilon^{-4})$ for just finding stationary points (see Appendix B and Table 1).

*Remark* 4.3. Offline methods are often studied under $(\varepsilon, \varepsilon^{1/2})$-approximate local minima. In the online setting, Ge et al. [22] used $(\varepsilon, \varepsilon^{1/4})$-approximate local minima, thus giving $T = O\big(\frac{\mathsf{poly}(d)}{\varepsilon^4} + \frac{\mathsf{poly}(d)}{\delta^{16}}\big)$. In general, it is better to treat $\varepsilon$ and $\delta$ separately to be more general, but nevertheless, $(\varepsilon, \varepsilon^{1/C})$-approximate local minima are always better than $\varepsilon$-approximate stationary points.

## 4.1 Our Theorem

We propose a new method $\texttt{Natasha2}^{\mathsf{full}}$ which, very informally speaking, alternatively

- finds approximate stationary points of $f(x)$ using `Natasha1.5`, or
- finds negative curvature of the Hessian $\nabla^2 f(x)$, using Oja's online eigenvector algorithm.

[11] The notion "$\nabla^2 f(x) \succeq -\delta \mathbf{I}$" means all the eigenvalues of $\nabla^2 f(x)$ are above $-\delta$.

[12] These methods are based on the "variance reduction" technique to reduce the random noise of SGD. They have been criticized by practitioners for performing poorer than SGD on training neural networks, because the noise of SGD allows it to escape from saddle points. Variance-reduction based methods have less noise and thus cannot escape from saddle points.

A similar alternation process (but for the offline problem) was studied by [1, 13]. Following their notion, we redefine gradient complexity $T$ to be the number of stochastic gradient computations plus Hessian-vector products.

Let $\Delta_f$ be any upper bound on $f(x_0) - f(x^*)$ where $x_0$ is the starting point. In this section, to present the simplest results, we use the big-$O$ notion to hide dependency in $L$, $L_2$, $\Delta_f$, and $\mathcal{V}$. In Section 7, we shall add back such dependency for a more general description of the algorithm. Our main result can be stated as follows:

**Theorem 2** (informal). *Under (A1), (A2) and (A4), for any $\varepsilon \in (0,1)$ and $\delta \in (0, \varepsilon^{1/4})$, $\texttt{Natasha2}(f, y_0, \varepsilon, \delta)$ outputs a point $x^{\mathsf{out}}$ so that, with probability at least $2/3$:*

$$\|\nabla f(x^{\mathsf{out}})\| \leq \varepsilon \quad \text{and} \quad \nabla^2 f(x^{\mathsf{out}}) \succeq -\delta \mathbf{I} \ .$$

*Furthermore, its gradient complexity is $T = \widetilde{O}\big(\frac{1}{\delta^5} + \frac{1}{\delta\varepsilon^3}\big)$ .*[13]

*Remark* 4.4. If $\delta > \varepsilon^{1/4}$ we can replace it with $\delta = \varepsilon^{1/4}$. Therefore, $T = \widetilde{O}\big(\frac{1}{\delta^5} + \frac{1}{\delta\varepsilon^3} + \frac{1}{\varepsilon^{3.25}}\big)$.

*Remark* 4.5. The follow-up work [10] replaced Hessian-vector products in `Natasha2` with only stochastic gradient computations, turning `Natasha2` into a pure first-order method. That paper is built on ours and thus all the proofs of this paper are still needed.

**Corollary 4.6.** *$T = \widetilde{O}(\varepsilon^{-3.25})$ for finding $(\varepsilon, \varepsilon^{1/4})$-approximate local minima. This is better than $T = O(\varepsilon^{-10/3})$ of SCSG for finding only $\varepsilon$-approximate stationary points.*

**Corollary 4.7.** *$T = \widetilde{O}(\varepsilon^{-3.5})$ for finding $(\varepsilon, \varepsilon^{1/2})$-approximate local minima. This was not known before this work, and is matched by several follow-up works using different algorithms [5, 10, 46, 48].*

### 4.2 Our Intuition

It is known that the problem of finding $(\varepsilon, \delta)$-approximate local minima, at a high level, "reduces" to (repeatedly) finding $\varepsilon$-approximate stationary points for an $O(\delta)$-nonconvex function [1, 13]. Specifically, Carmon et al. [13] proposed the following procedure. In every iteration at point $y_k$, detect whether the minimum eigenvalue of $\nabla^2 f(y_k)$ is below $-\delta$:

- if yes, find the minimum eigenvector of $\nabla^2 f(y_k)$ approximately and move in this direction.
- if no, let $F^k(x) \stackrel{\text{def}}{=} f(x) + L\big(\max\big\{0, \|x - y_k\| - \frac{\delta}{L_2}\big\}\big)^2$, which can be proven as $5L$-smooth and $3\delta$-nonconvex; then find an $\varepsilon$-approximate stationary point of $F^k(x)$ to move there. Intuitively, $F^k(x)$ penalizes us from moving out of the "safe zone" of $\big\{x \colon \|x - y_k\| \leq \frac{\delta}{L_2}\big\}$.

Previously, it was thought necessary to achieve high accuracy for both tasks above. This is why researchers have only been able to design offline methods: in particular, the shift-and-invert method [21] was applied to find the minimum eigenvector, and `repeatSVRG` was applied to find a stationary point of $F^k(x)$.[14]

In this paper, we apply efficient *online* algorithms for the two tasks: namely, Oja's algorithm (see Section 5.1) for finding minimum eigenvectors, and our new `Natasha1.5` algorithm (see Section 3.2)

[13] Throughout this paper, we use the $\widetilde{O}$ notion to hide at most one logarithmic factor in all the parameters (namely, $n, d, L, L_2, \mathcal{V}, 1/\varepsilon, 1/\delta$).

[14] The performance of `repeatSVRG` was summarized in Table 2 and Figure 4(a). `repeatSVRG` is an offline algorithm, and finds an $\varepsilon$-approximate stationary point for a function $f(x)$ that is $\sigma$-nonconvex. It is divided into stages. In each stage $t$, it considers a modified function $f_t(x) \stackrel{\text{def}}{=} f(x) + \sigma\|x - x_t\|^2$, and then apply the accelerated `SVRG` method to minimize $f_t(x)$. Then, it moves to $x_{t+1}$ which is a sufficiently accurate minimizer of $f_t(x)$.

**Algorithm 2** `Natasha2`$(f, y_0, \varepsilon, \delta)$

---

**Input:** function $f(x) = \frac{1}{n}\sum_{i=1}^n f_i(x)$, starting vector $y_0$, target accuracy $\varepsilon > 0$ and $\delta > 0$.

1: **if** $\frac{\varepsilon^{1/3}}{\delta} \geq 1$ **then** $\widetilde{L} = \widetilde{\sigma} \leftarrow \Theta(\frac{\varepsilon^{1/3}}{\delta}) \geq 1$; $\diamond$ *the boundary case for large* $L_2$
2: **else** $\widetilde{L} \leftarrow 1$ and $\widetilde{\sigma} \leftarrow \Theta\big(\frac{\varepsilon}{\delta^3}\big) \in [\delta, 1]$.
3: $X \leftarrow []$;
4: **for** $k \leftarrow 0$ **to** $\infty$ **do**
5: Apply Oja's algorithm to find minEV $v$ of $\nabla^2 f(y_k)$ for $\widetilde{\Theta}(\frac{1}{\delta^2})$ iterations $\diamond$ *see Lemma 5.3*
6: **if** $v \in \mathbb{R}^d$ is found s.t. $v^\top \nabla^2 f(y_k) v \leq -\frac{\delta}{2}$ **then**
7: $y_{k+1} \leftarrow y_k \pm \frac{\delta}{L_2} v$ where the sign is random.
8: **else** $\diamond$ *it satisfies* $\nabla^2 f(y_k) \succeq -\delta\mathbf{I}$
9: $F^k(x) \stackrel{\text{def}}{=} f(x) + L(\max\{0, \|x - y_k\| - \frac{\delta}{L_2}\})^2$.
10: run `Natasha1.5`$\big(F^k, y_k, \Theta(\varepsilon^{-2}), 1, \Theta(\varepsilon\delta)\big)$. $\diamond$ $F^k(\cdot)$ *is* $\widetilde{L}$*-smooth and* $\widetilde{\sigma}$*-nonconvex*
11: let $\widehat{\mathsf{y}}_k, y_{k+1}$ be the vector $\widehat{\mathsf{y}}$ and $\widehat{\mathsf{x}}$ when Line 13 is reached in `Natasha1.5`.
12: $X \leftarrow [X, (y_k, \widehat{\mathsf{y}}_k)]$;
13: break the for loop if have performed $\Theta\big(\frac{1}{\delta\varepsilon}\big)$ first-order steps.
14: **end if**
15: **end for**
16: $(y, \widehat{\mathsf{y}}) \leftarrow$ a random pair in $X$. $\diamond$ *in practice, simply output* $\widehat{\mathsf{y}}_k$
17: define convex function $g(x) \stackrel{\text{def}}{=} f(x) + L(\max\{0, \|x - y\| - \frac{\delta}{L_2}\})^2 + \widetilde{\sigma}\|x - \widehat{\mathsf{y}}\|^2$.
18: use SGD to minimize $g(x)$ for $\widetilde{\Theta}(\frac{1}{\varepsilon^2})$ steps and output $x^{\mathsf{out}}$.

---

for finding stationary points. More specifically, for Oja's, we only decide if there is an eigenvalue below threshold $-\delta/2$, or conclude that the Hessian has all eigenvalues above $-\delta$. This can be done in an online fashion using $O(\delta^{-2})$ Hessian-vector products (with high probability) using Oja's algorithm. For `Natasha1.5`, we only apply it for a single epoch of length $B = \Theta(\varepsilon^{-2})$. Conceptually, this shall make the above procedure online and run in a complexity independent of $n$.

*Unfortunately,* technical issues arise in this "wishful thinking."

Most notably, the above process finishes only if `Natasha1.5` finds an approximate stationary point $x$ of $F^k(x)$ that is *also* inside the safe zone $\big\{x \colon \|x - y_k\| \leq \frac{\delta}{L_2}\big\}$. This is because $F^k(x) = f(x)$ inside the safe zone and therefore $\|\nabla F^k(x)\| \leq \varepsilon$ also implies $\|\nabla f(x)\| \leq 2\varepsilon$.

What can we do if we move out of the safe zone? To tackle this case, we show an additional property of `Natasha1.5` (see Lemma 6.5). That is, the amount of objective decrease —i.e., $f(y_k) - f(x)$ if $x$ moves out of the safe zone— must be proportional to the distance $\|x - y_k\|^2$ we travel in space. Therefore, if $x$ moves out of the safe zone, then we can decrease sufficiently the objective. This is also a good case. This summarizes some high-level technical ingredient of `Natasha2`.

We formally state `Natasha2` in Algorithm 2, and it uses the big-$O$ notion to hide dependency in $L$, $L_2$, $\mathcal{V}$ and $\Delta_f$. The more general code to take care of all the parameters can be found in Algorithm 5 of Section 7. We point out that Remark 3.4 and Remark 3.5 still apply here (regarding why we need to randomly select $\widehat{\mathsf{y}}$ or to do pruning in Line 17 of `Natasha2`).

Finally, we stress that although we borrowed the construction of $f(x) + L\big(\max\big\{0, \|x - y_k\| - \frac{\delta}{L_2}\big\}\big)^2$ from the offline algorithm of Carmon et al. [13], our `Natasha2` algorithm and analysis are different from them in all other aspects.

# 5 Auxiliary Lemmas

We show a few auxiliary results that shall be used later in the analysis of $\texttt{Natasha1.5}^{\mathsf{full}}$ and $\texttt{Natasha2}^{\mathsf{full}}$.

- In Section 5.1, we revisit Oja's algorithm which is an online method for finding eigenvectors.
- In Section 5.2, we present a new sufficient condition for finding stationary points.
- In Section 5.3, we recall a few results for SGD on convex functions.

## 5.1 Oja's Algorithm

Let $\mathcal{D}$ be a distribution over $d \times d$ symmetric matrices whose eigenvalues are between 0 and 1, and denote by $\mathbf{B} \stackrel{\mathrm{def}}{=} \mathbb{E}_{\mathbf{A}\sim\mathcal{D}}[\mathbf{A}]$ its mean. Let $\mathbf{A}_1, \dots, \mathbf{A}_T$ be $T$ copies of i.i.d. samples generated from $\mathcal{D}$. Oja's algorithm begins with a random unit-norm Gaussian vector $w_1 \in \mathbb{R}^d$. At each iteration $k \in 2, \dots, T$, Oja's algorithm computes $w_k = \frac{(\mathbf{I}+\eta\mathbf{A}_{k-1})w_{k-1}}{C}$ where $C > 0$ is the normalization constant such that $\|w_k\| = 1$. Allen-Zhu and Li [9] showed (see its last section) that [15]

**Theorem 5.1.** *For every $p \in (0,1)$, choosing $\eta = \Theta(\sqrt{p/T})$, we have with prob. $\geq 1-p$:*

$$\sum_{k=1}^T w_k^\top \mathbf{B} w_k \geq T \cdot \lambda_{\max}(\mathbf{B}) - O\big(\tfrac{\sqrt{T}}{\sqrt{p}} \cdot \log(d/p)\big) \enspace .$$

*Remark* 5.2. The above result is asymptotically optimal even in terms of sampling complexity [9].

**Approximating MinEV of Hessian.** Suppose $f(x) = \frac{1}{n}\sum_{i=1}^n f_i(x)$ where each $f_i(x)$ is twice-differentiable and $L$-smooth. We can denote by $\mathcal{D}$ the distribution where each $\frac{L\cdot\mathbf{I} - \nabla^2 f_i(x)}{2L} \in \mathbb{R}^{d\times d}$ is generated with probability $\frac{1}{n}$, and then use Oja's algorithm to compute the minimum eigenvalue of $\nabla^2 f(x)$. Note that each time when computing $(\mathbf{I}+\eta\mathbf{A}_{k-1})w_{k-1}$, it suffices to compute Hessian-vector product (i.e., $\nabla^2 f_i(x)\cdot w_{k-1}$) once. The following corollary is simple to prove:

**Lemma 5.3.** *There exists absolute constant $C > 1$ such that for any $x \in \mathbb{R}^d$, $T \geq 1$, $p \in (0,1)$:*

- *if we run Oja's algorithm once for $T$ iterations, with $\eta = \Theta(\sqrt{T})$, we can find unit vector $y$ such that, with at with probability at least $4/5$,*

$$y^\top \nabla^2 f(x) y \leq \lambda_{\min}(\nabla^2 f(x)) + C \cdot \tfrac{L\log(d)}{\sqrt{T}} \enspace .$$

- *if we run Oja's algorithm $O(\log(1/p))$ times each for $T$ iterations, then w.p. $\geq 1-p$, we can*

$$\text{either conclude} \qquad \lambda_{\min}(\nabla^2 f(x)) \geq -C \cdot \frac{L\log(d/p)}{\sqrt{T}} \enspace ,$$

$$\text{or find } y \in \mathbb{R}^d\text{:} \qquad y^\top \nabla^2 f(x) y \leq -\frac{C}{2} \cdot \frac{L\log(d/p)}{\sqrt{T}} \enspace .$$

*The total number of Hessian-vector products is at most $O(T\log(1/p))$.*

*Remark* 5.4. We refer to the computation of $\nabla^2 f_i(x)\cdot v$ for $i \in [n]$ and $v \in \mathbb{R}^d$ as a Hessian-vector product. Therefore, computing $\nabla^2 f(x)\cdot v$ counts as $n$ times of Hessian-vector products.

In a follow-up work, Allen-Zhu and Li [10] designed a minor variant of Oja's algorithm which achieves the same guarantee as Lemma 5.3 but using only stochastic gradient computations (without Hessian-vector products). The main idea is to use $\frac{\nabla f_i(x+qw)-\nabla f_i(x)}{q}$ to replace the use of $\nabla^2 f_i(x)\cdot w$ for some small constant $q > 0$.

[15] The original one-paged proof from [9] only showed Theorem 5.1 where the left hand side is $\sum_{k=1}^T w_k^\top \mathbf{A}_k w_k$. However, by Azuma's inequality, we have $\sum_{k=1}^T w_k^\top \mathbf{B} w_k \geq \sum_{k=1}^T w_k^\top \mathbf{A}_k w_k - O(\sqrt{T\log(1/p)})$ with probability $\geq 1-p$.

## 5.2 First-Order Stopping Criterion

We present a sufficient condition for finding approximate stationary points for

$$F(x) = \psi(x) + f(x) \ , \tag{5.1}$$

where $\psi(x)$ is proper convex, $f(x)$ is $\sigma$-*nonconvex* but $L$-smooth. For any $\widehat{\mathsf{x}} \in \mathbb{R}^d$, if we define

$$G(x) \stackrel{\mathrm{def}}{=} \psi(x) + g(x) \stackrel{\mathrm{def}}{=} \psi(x) + \big(f(x) + \sigma\|x - \widehat{\mathsf{x}}\|^2\big) \ ,$$

then $g(x)$ becomes $\sigma$-strongly convex, and thus we can use convex optimization to minimize $G(x)$.

The following lemma says that, if we find an approximate stationary point $x$ of $G(x)$, then it is also an approximate stationary point of $F(x)$ up to an additive error $O(\sigma^2\|\widehat{\mathsf{x}} - x^*\|^2)$, where $x^*$ is the exact minimizer of $G(x)$.

**Lemma 5.5.** *Let $x^*$ be the unique minimizer of $G(y)$, and $x$ be an arbitrary vector in the domain of $\{x \in \mathbb{R}^d \colon \psi(x) < +\infty\}$. Then, for every $\eta \in \big(0, \frac{1}{L+2\sigma}\big]$, we have*

$$\|\mathcal{G}_{F,\eta}(x)\|^2 + \sigma^2\|x - \widehat{\mathsf{x}}\|^2 \le O\big(\sigma^2\|x^* - \widehat{\mathsf{x}}\|^2 + \|\mathcal{G}_{G,\eta}(x)\|^2\big) \ .$$

*Remark* 5.6. When $\psi(x) \equiv 0$ and $x = x^*$, Lemma 5.5 is trivial: $\|\mathcal{G}_{F,\eta}(x)\| = \|\nabla F(x)\| = \|\nabla G(x) - 2\sigma(x - \widehat{\mathsf{x}})\| = 2\sigma\|x - \widehat{\mathsf{x}}\|$. The main difficulty arises in order to deal with $\psi(x) \neq 0$ and $x \neq x^*$.

Let us compare Lemma 5.5 to its close variant shown in the work of `Natasha1` [3]. In [3], the author proved a similar result as Lemma 5.5, with $\|\mathcal{G}_{G,\eta}(x)\|^2$ replaced by $\frac{G(x)-G(x^*)}{\eta^2\sigma}$. The result in [3] is weaker, because even if $\psi(x) = 0$ and even if $\eta = 1/(L + 2\sigma)$, we have

$$\|\mathcal{G}_{G,\eta}(x)\|^2 = \|\nabla G(x)\|^2 \le L(G(x) - G(x^*)) \ll \tfrac{1}{\eta^2\sigma}(G(x) - G(x^*)) \ .$$

If using this weaker version, our convergence rate shall become worsened.[16]

## 5.3 Proximal SGD for Convex Optimization

We revisit stochastic gradient descent (SGD) on minimizing a *convex* stochastic objective

$$F(x) = \psi(x) + f(x) \stackrel{\mathrm{def}}{=} \psi(x) + \tfrac{1}{n}\textstyle\sum_{i\in[n]} f_i(x) \ , \tag{5.2}$$

where

1. $\psi(x)$ is proper convex,
2. each $f_i(x)$ is differentiable, $f(x)$ is convex and $L$-smooth,
3. $F(x)$ is $\sigma$-strongly convex for some $\sigma \in [0, L]$, and
4. the stochastic gradients $\nabla f_i(x)$ have a bounded variance (over the domain of $\psi(\cdot)$), that is

$$\forall x \in \{y \in \mathbb{R}^d \,|\, \psi(y) < +\infty\} \colon \quad \mathbb{E}_{i\in_R[n]}\|\nabla f(x) - \nabla f_i(x)\|^2 \le \mathcal{V} \ .$$

Recall that SGD repeatedly performs *proximal updates* of the form

$$x_{t+1} = \arg\min_{y\in\mathbb{R}^d}\{\psi(y) + \tfrac{1}{2\alpha}\|y - x_t\|^2 + \langle\nabla f_i(x_t), y\rangle\} \ ,$$

where $\alpha > 0$ is some learning rate, and $i \in_R [n]$ per iteration. Note that if $\psi(y) \equiv 0$ then $x_{t+1} = x_t - \alpha\nabla f_i(x_t)$. Let gradient complexity $T$ be the number of computations of $\nabla f_i(x)$.

The next theorem is due to Allen-Zhu [5, Theorem 3]:

[16] Using Lemma 5.5, to find $\varepsilon$-approximate stationary points of $F(x)$, we wish to find a point $x$ satisfying $\|\mathcal{G}_{G,\eta}(x)\|^2 \le \varepsilon^2$. The convergence rate for SGD to achieve this goal is $O(\frac{1}{\varepsilon^2})$, see Theorem 5.7. In contrast, if using [3], one needs to find $G(x) - G(x^*) \le O(\sigma\varepsilon^2)$. The convergence rate to achieve this goal is $O(\frac{1}{\sigma^2\varepsilon^2})$. This worse dependency on $\sigma$ shall slow down the performance of our proposed methods `Natasha1.5` and `Natasha2`.

**Theorem 5.7** ([5]). *Let $x^* \in \arg\min_x\{F(x)\}$ and $C \in (0,1]$ be any absolute constant. To solve Problem (5.2) given a starting vector $x_0 \in \mathbb{R}^d$, there is an SGD variant* $\mathtt{SGD3}^{\mathsf{sc}}$ *which, for every $T \geq \frac{L}{\sigma}\log\frac{L}{\sigma}$,* $\mathtt{SGD3}^{\mathsf{sc}}(F, x_0, \sigma, L, T)$ *computes $T$ stochastic gradients and outputs $\overline{x}$ with*

$$\mathbb{E}[\|\mathcal{G}_{F,\eta}(\overline{x})\|] \leq O\Big(\frac{\sqrt{\mathcal{V}}\cdot\log^{3/2}\frac{L}{\sigma}}{\sqrt{T}}\Big) + \big(1-\frac{\sigma}{L}\big)^{\Omega(T/\log(L/\sigma))}\sigma\|x_0 - x^*\| \quad \textit{where } \eta = C/L \ .$$

In other words, to find a point with $\|\mathcal{G}_{F,\eta}(\overline{x})\| \leq \varepsilon$, the method $\mathtt{SGD3}^{\mathsf{sc}}$ needs $T \propto \widetilde{O}(\varepsilon^{-2})$ stochastic gradient computations. In contrast, the naive SGD gives only $T \propto \widetilde{O}(\sigma^{-1}\varepsilon^{-2})$. (See the discussions in [5] and the references therein.) We use this better rate $T \propto \widetilde{O}(\varepsilon^{-2})$ in order to tighten our final complexities of `Natasha1.5` and `Natasha2`.

# 6 Natasha 1.5: Finding Stationary Points

In this section, we study the problem finding approximate stationary points for

$$F(x) \stackrel{\text{def}}{=} \psi(x) + f(x) \stackrel{\text{def}}{=} \psi(x) + \tfrac{1}{n}\textstyle\sum_{i=1}^n f_i(x) \ , \tag{6.1}$$

where

1. $\psi(\cdot)$ is proper convex,
2. each $f_i(x)$ is possibly nonconvex but $L$-smooth,
3. the average $f(x)$ is $\sigma$-nonconvex for $\sigma \in (0, L]$,[17] and
4. the stochastic gradients $\nabla f_i(x)$ have a bounded variance (over the domain of $\psi(\cdot)$), that is

$$\forall x \in \{y \in \mathbb{R}^d \,|\, \psi(y) < +\infty\}: \quad \mathbb{E}_{i\in_R[n]}\|\nabla f(x) - \nabla f_i(x)\|^2 \leq \mathcal{V} \ .$$

Throughout this section, we define $T$, the gradient complexity, as the number of computations of $\nabla f_i(x)$. For simplicity, we explain the intuition in the special case when $\psi(x) \equiv 0$.

**Algorithm.** Our full pseudocode $\mathtt{Natasha1.5}^{\mathsf{full}}$ is given in Algorithm 3. It consists of $T'$ full epochs $k = 1, \dots, T'$. At the beginning of each full epoch, we compute $\nabla f_S(\widetilde{\mathsf{x}}) \stackrel{\text{def}}{=} \frac{1}{|S|}\sum_{i\in S}\nabla f_i(x)$ for a random subset $S \subseteq [n]$ of cardinality $B$, where $\widetilde{\mathsf{x}}$ is the current snapshot point.

Each epoch $k$ is is further divided into $p$ sub-epochs $s = 0, 1, \dots, p-1$, each of length $m = B/p$. In each sub-epoch $s$, we start with a point $x_0 = \widehat{\mathsf{x}}$, and conceptually apply `SVRG` but replacing $f(x)$ with its regularized version $f^s(x) \stackrel{\text{def}}{=} f(x) + \sigma\|x - \widehat{\mathsf{x}}\|^2$. In other words, we

- compute gradient estimator $\widetilde{\nabla} = \nabla f_S(\widetilde{\mathsf{x}}) + \nabla f_i(x_t) - \nabla f_i(x_t) + 2\sigma(x_t - \widehat{\mathsf{x}})$, and
- perform update $x_{t+1} = \arg\min_y \big\{\psi(y) + \langle\widetilde{\nabla}, y\rangle + \frac{1}{2\alpha}\|y - x_t\|^2\big\}$ with learning rate $\alpha$.

Finally, when the sub-epoch is over, we define $\widehat{\mathsf{x}}$ to be a random one from $\{x_0, \dots, x_{m-1}\}$; when a full epoch is over, we define $\widetilde{\mathsf{x}}$ to be the last $\widehat{\mathsf{x}}$.

In the end, we output two points for later use, $\widehat{\mathsf{y}}$ is a random $\widehat{\mathsf{x}}$ among all the full epochs and sub-epochs, and $y^+$ is the last $\widehat{\mathsf{x}}$. Very informally speaking, $\|\nabla f(\widehat{\mathsf{y}})\|$ is roughly upper bounded by $f(x^\varnothing) - f(y^+)$; in other words, $\widehat{\mathsf{y}}$ is a point that gives small gradient, but $y^+$ is a point that ensures objective decrease

We analyze the behavior of $\mathtt{Natasha1.5}^{\mathsf{full}}$ for one full epoch in Section 6.1 and then telescope it for all epochs in Section 6.2.

[17] We assume $\sigma \leq L$ without loss of generality, because any $L$-smooth function is also $L$-nonconvex.

**Algorithm 3** $\texttt{Natasha1.5}^{\mathsf{full}}(F, x^{\varnothing}, B, p, T', \alpha)$

**Input:** function $F(\cdot)$ satisfying Problem (6.1), starting vector $x^{\varnothing}$, epoch length $B \in [n]$, sub-epoch count $p \in [B]$, epoch count $T' \geq 1$, learning rate $\alpha > 0$. $\diamond$ *p should be* $\Theta((\sigma^2 B/L^2)^{1/3})$

**Output:** two vectors $\widehat{\mathsf{y}}$ and $y^+$.

```
x̂ ← x^∅; m ← B/p; X ← [];
for k ← 1 to T' do                                                   ◇ T' epochs each of length B
    x̃ ← x̂; μ ← (1/B) Σ_{i∈S} ∇f_i(x̃) where S is a uniform random subset of [n] with |S| = B;
    for s ← 0 to p − 1 do                                            ◇ p sub-epochs each of length m
        x_0 ← x̂; X ← [X, x̂];
        for t ← 0 to m − 1 do                                        ◇ m iterations in each sub-epoch
            i ← a random index from [n].
            ∇̃ ← ∇f_i(x_t) − ∇f_i(x̃) + μ + 2σ(x_t − x̂)
            x_{t+1} = arg min_{y∈R^d} {ψ(y) + (1/2α)‖y − x_t‖² + ⟨∇̃, y⟩}
        end for
        x̂ ← a random choice from {x_0, x_1, ..., x_{m−1}};            ◇ in practice, choose the average
    end for
end for
ŷ ← a random vector in X and y⁺ ← x̂.                                 ◇ in practice, choose the last
return (ŷ, y⁺).
```

## 6.1 Natasha 1.5: Analysis for One Epoch

**Notations.** When focusing on a single full epoch (with $k$ being fixed), we introduce the following notations for analysis purpose only.

- Let $\widehat{\mathsf{x}}^s$ be the vector $\widehat{\mathsf{x}}$ at the beginning of sub-epoch $s$.
- Let $x_t^s$ be the vector $x_t$ in sub-epoch $s$.
- Let $i_t^s$ be the index $i \in [n]$ in sub-epoch $s$ at iteration $t$.
- Let $f^s(x) \stackrel{\text{def}}{=} f(x) + \sigma\|x - \widehat{\mathsf{x}}^s\|^2$, $F^s(x) \stackrel{\text{def}}{=} F(x) + \sigma\|x - \widehat{\mathsf{x}}^s\|^2$, and $x_*^s \stackrel{\text{def}}{=} \arg\min_x\{F^s(x)\}$.
- Let $\widetilde{\nabla} f^s(x_t^s) \stackrel{\text{def}}{=} \nabla f_i(x_t^s) - \nabla f_i(\widetilde{\mathbf{x}}) + \nabla f_S(\widetilde{\mathbf{x}}) + 2\sigma(x_t - \widehat{\mathsf{x}})$ where $i = i_t^s$.
- Let $\widetilde{\nabla} f(x_t^s) \stackrel{\text{def}}{=} \nabla f_i(x_t^s) - \nabla f_i(\widetilde{\mathbf{x}}) + \nabla f_S(\widetilde{\mathbf{x}})$ where $i = i_t^s$.
- Let $\mathbf{e} \stackrel{\text{def}}{=} \nabla f_S(\widetilde{\mathbf{x}}) - \nabla f(\widetilde{\mathbf{x}})$.

We obviously have that $f^s(x)$ and $F^s(x)$ are $\sigma$-strongly convex, and $f^s(x)$ is $(L + 2\sigma)$-smooth.

The following lemma gives an upper bound on the variance of the gradient estimator $\widetilde{\nabla} f^s(x_t^s)$. The only difference to $\texttt{Natasha1}$ [3] is the additional term $\|\mathbf{e}\|^2$.

**Lemma 6.1.** *We have* $\mathbb{E}_{i_t^s}\big[\|\widetilde{\nabla} f^s(x_t^s) - \nabla f^s(x_t^s)\|^2\big] \leq pL^2\|x_t^s - \widehat{\mathsf{x}}^s\|^2 + pL^2\sum_{k=0}^{s-1}\|\widehat{\mathsf{x}}^k - \widehat{\mathsf{x}}^{k+1}\|^2 + \|\mathbf{e}\|^2$ .

The following simple claim bounds $\|\mathbf{e}\|^2$.

**Claim 6.2.** *If $S$ is a uniform random subset of $[n]$ with cardinality $|S| = B$, then* $\mathbb{E}_S[\|\mathbf{e}\|^2] \leq \frac{\mathcal{V}}{B}$.

*Proof of Claim 6.2.* We first recall that if $v_1, \dots, v_n \in \mathbb{R}^d$ satisfy $\sum_{i=1}^n v_i = \vec{0}$, and $S$ is a non-empty, uniform random subset of $[n]$. Then

$$\mathbb{E}\Big[\big\|\tfrac{1}{|S|}\textstyle\sum_{i\in S} v_i\big\|^2\Big] = \tfrac{n-|S|}{(n-1)|S|} \cdot \tfrac{1}{n}\textstyle\sum_{i\in[n]}\|v_i\|^2 \leq \tfrac{\mathbb{I}[|S|<n]}{|S|} \cdot \tfrac{1}{n}\textstyle\sum_{i\in[n]}\|v_i\|^2 \ .$$

Letting $v_i = \nabla f_i(\widetilde{\mathbf{x}}) - \nabla f(\widetilde{\mathbf{x}})$, we have

$$\mathbb{E}[\|\mathbf{e}\|^2] = \mathbb{E}\Big[\|\frac{1}{|S|}\sum_{i\in S} v_i\|^2\Big] \leq \frac{\mathbb{I}[|S|<n]}{|S|}\cdot\frac{1}{n}\sum_{i\in[n]}\|\nabla f_i(\widetilde{\mathbf{x}}) - \nabla f(\widetilde{\mathbf{x}})\|^2 \leq \frac{\mathcal{V}}{B} \ . \qquad \square$$

The following lemma is our main contribution for the base method $\texttt{Natasha1.5}^{\mathsf{full}}$. It is analogous to the main lemma of Natasha1 [3]; however, we have to apply additional tricks to handle the fact that $\widetilde{\nabla} f^s(x)$ is a *biased* estimator of $\nabla f^s(x)$. (Recall that $\mathbb{E}_{i_t^s}[\widetilde{\nabla} f^s(x_t^s)] = \nabla f^s(x_t^s) + \mathbf{e}$.)

*Remark* 6.3. The proof of Lemma 6.4 only relies on mirror descent. This is different from the gradient-descent analysis of SCSG [30], and thus very different from how the proof of SCSG handles this additional bias $\mathbf{e}$. We believe this is the key for achieving our result on $\texttt{Natasha1.5}^{\mathsf{full}}$.

**Lemma 6.4.** *As long as* $\alpha \leq \frac{1}{2L+4\sigma}$, *letting* $x_*^s = \arg\min_x\{F(x) + \sigma\|x - \widehat{\mathsf{x}}^s\|^2\}$, *we have*

$$\mathbb{E}\Big[\big(F^s(\widehat{\mathsf{x}}^{s+1}) - F^s(x_*^s)\big)\Big] \leq \mathbb{E}\Big[\frac{F^s(\widehat{\mathsf{x}}^s) - F^s(x_*^s)}{\sigma\alpha m/4} + 2\alpha p L^2\Big(\sum_{k=0}^{s}\|\widehat{\mathsf{x}}^k - \widehat{\mathsf{x}}^{k+1}\|^2\Big)\Big] + \frac{3}{\sigma}\|\mathbf{e}\|^2 \ .$$

One can telescope Lemma 6.4 for an entire epoch and arrive at the following lemma:

**Lemma 6.5.** *If* $\alpha \leq \frac{1}{2L+4\sigma}$, $\alpha \geq \frac{8}{\sigma m}$ *and* $\alpha \leq \frac{\sigma}{4p^2L^2}$, *we have*

$$\sum_{s=0}^{p-1}\mathbb{E}\Big[\sigma\|\widehat{\mathsf{x}}^s - \widehat{\mathsf{x}}^{s+1}\|^2 + \frac{\sigma}{2}\|\widehat{\mathsf{x}}^s - x_*^s\|^2\Big] \leq 2\mathbb{E}\Big[F(\widehat{\mathsf{x}}^0) - F(\widehat{\mathsf{x}}^p)\Big] + \frac{3p\mathcal{V}}{\sigma B} \ ,$$

*where recall* $x_*^s \stackrel{\text{def}}{=} \arg\min_x\{F(x) + \sigma\|x - \widehat{\mathsf{x}}^s\|^2\}$.

## 6.2 Natasha 1.5: Final Theorem

As we shall see in the next section, the design of $\texttt{Natasha2}^{\mathsf{full}}$ for finding approximate local minima requires to run $\texttt{Natasha1.5}^{\mathsf{full}}$ only for *one* full epoch, that is, $T' = 1$. However, for the purpose of achieving good stationary points and proving Theorem 1, we need to run $\texttt{Natasha1.5}^{\mathsf{full}}$ for $T' \geq 1$ and then apply SGD for pruning in the end.

Specifically, as summarized in $\texttt{Natasha1.5}^{\mathsf{prune}}$, we specify parameters $B$, $p$, and $\alpha$ appropriately and call $\texttt{Natasha1.5}^{\mathsf{full}}$. Then, we perform an additional SGD starting from $\widehat{\mathsf{y}}$ and output $x^{\mathsf{out}}$.

---

**Algorithm 4** $\texttt{Natasha1.5}^{\mathsf{prune}}(F, x^{\varnothing}, \varepsilon)$

---

**Input:** function $F(\cdot)$ satisfying Problem (6.1), starting vector $x^{\varnothing}$, either gradient complexity $T \geq 1$ or target accuracy $\varepsilon > 0$.

1: $B \leftarrow \Theta(\frac{\mathcal{V}}{\varepsilon^2})$; $p \leftarrow (\frac{\sigma^2}{48L^2}B)^{1/3}$; $\alpha \leftarrow \Theta(\frac{\sigma}{p^2L^2})$. ◇ *$p \in [1, B]$ under assumption $L \geq \sigma \geq \Omega(\frac{\varepsilon L}{\mathcal{V}^{1/2}})$*

2: $(\widehat{\mathsf{y}}, y^+) \leftarrow \texttt{Natasha1.5}^{\mathsf{full}}(F, x^{\varnothing}, B, p, T/B, \alpha)$ for $T = \Theta\big(\frac{(L^2\sigma)^{1/3}\Delta_F\cdot\mathcal{V}^{2/3}}{\varepsilon^{10/3}}\big)$;
◇ *$\Delta_F$ is an upper bound on $F(x^{\varnothing}) - \min_x\{F(x)\}$.*

3: define convex function $G(x) \stackrel{\text{def}}{=} F(x) + \sigma\|x - \widehat{\mathsf{y}}\|^2$.

4: $x^{\mathsf{out}} \leftarrow \texttt{SGD3}^{\mathsf{sc}}(G, \widehat{\mathsf{y}}, \sigma, O(L), T_{\mathsf{sgd}})$ for $T_{\mathsf{sgd}} = \Theta\big(\frac{\mathcal{V}}{\varepsilon^2}\log^3\frac{L}{\sigma}\big)$.

5: **return** $x^{\mathsf{out}}$.

---

We are now ready to state and prove our main convergence theorem for $\texttt{Natasha1.5}^{\mathsf{full}}$:

**Theorem 1.** *Consider Problem (6.1) with a starting vector $x^{\varnothing}$. Let $\eta = \frac{C}{L}$ where $C \in (0,1]$ is any absolute constant, and $\Delta_F$ be an upper bound on $F(x^{\varnothing}) - \min_x\{F(x)\}$. If $\varepsilon > 0$ and $L \geq \sigma \geq \Omega\big(\frac{\varepsilon L}{\mathcal{V}^{1/2}}\big)$, then* $\mathtt{Natasha1.5}^{\mathsf{prune}}(F, x^{\varnothing}, \varepsilon)$ *computes an output $x^{\mathsf{out}}$ with $\mathbb{E}[\|\mathcal{G}_{F,\eta}(x^{\mathsf{out}})\|] \leq \varepsilon$ in gradient complexity*

$$T = O\Big(\frac{\mathcal{V}}{\varepsilon^2}\log^3\frac{L}{\sigma} + \frac{(L^2\sigma)^{1/3}\Delta_F\cdot\mathcal{V}^{2/3}}{\varepsilon^{10/3}}\Big) \enspace .$$

Since we can always replace $\sigma$ with $\frac{\varepsilon L}{\mathcal{V}^{1/2}}$ if it is too small, we have the following corollary:

**Corollary 6.6.** *Treating $L$, $\Delta_F$, and $\mathcal{V}$ as constants, we have $T = O\big(\frac{1}{\varepsilon^3} + \frac{\sigma^{1/3}}{\varepsilon^{10/3}}\big)$.*

*Remark* 6.7. If $\varepsilon$ is not known before the execution, one can similarly write $\mathtt{Natasha1.5}^{\mathsf{prune}}(F, x^{\varnothing}, T)$ in terms of an input parameter $T$. This requires setting $B = \Theta\big(\frac{\mathcal{V}^3 T^3}{\Delta_F^3 L^2\sigma}\big)^{1/5}$ rounded to the nearest integer between 1 and $T$. One can prove that, if $T \geq \Omega\big(\frac{L}{\sigma}\log^2\frac{L}{\sigma} + \frac{L^4\Delta_F}{\sigma^3\mathcal{V}}\big)$, then in gradient complexity $T$, $\mathtt{Natasha1.5}^{\mathsf{prune}}(F, x^{\varnothing}, T)$ outputs a point $x^{\mathsf{out}}$ with

$$\mathbb{E}[\|\mathcal{G}_{F,\eta}(x^{\mathsf{out}})\|] \leq O\Big(\frac{\sqrt{\mathcal{V}}\log^{3/2}\frac{L}{\sigma} + L^{1/3}\sigma^{1/6}\Delta_F^{1/2}}{\sqrt{T}} + \frac{\sigma^{1/10}L^{1/5}\Delta_F^{3/10}\mathcal{V}^{1/5}}{T^{3/10}}\Big) \enspace .$$

*Proof of Theorem 1.* Recall we choose

$$p \stackrel{\text{def}}{=} \big(\tfrac{\sigma^2}{48L^2}B\big)^{1/3} \in [1, B] \quad , \quad m \stackrel{\text{def}}{=} B/p \quad \text{and} \quad \alpha \stackrel{\text{def}}{=} \tfrac{8}{\sigma m} = \tfrac{\sigma}{6p^2L^2} \leq \tfrac{\sigma}{6L^2} \leq \tfrac{1}{2L+4\sigma} \enspace .$$

These parameters satisfy the prerequisite of Lemma 6.5. We denote by $T' = T/B$.

If we telescope Lemma 6.5 for the entire algorithm (which has $T'$ full epochs), and use the fact that $\widehat{\mathsf{x}}^p$ of the previous epoch equals $\widehat{\mathsf{x}}^0$ of the next epoch, we conclude that if we choose $\widehat{\mathsf{y}}$ to be $\widehat{\mathsf{x}}_s$ for a random epoch and a random subepoch $s$, and $y^+ = \widehat{\mathsf{x}}^p$ of the last epoch, we have

$$\mathbb{E}[\sigma\|\widehat{\mathsf{y}} - \widehat{\mathsf{y}}^*\|^2] \leq \frac{2}{pT'}(F(x^{\varnothing}) - \mathbb{E}[F(y^+)]) + \frac{3\mathcal{V}}{\sigma B}$$

where recall $\widehat{\mathsf{y}}^* = \arg\min_y\{F(y) + \sigma\|x - \widehat{\mathsf{y}}\|^2\}$. By the choice $T' = T/B$, the choice of $p$, and $F(y^+) \geq \min_x\{F(x)\}$, we have

$$\mathbb{E}[\sigma\|\widehat{\mathsf{y}} - \widehat{\mathsf{y}}^*\|^2] \leq O\Big(\frac{L^{2/3}B^{2/3}}{\sigma^{2/3}T}(F(x^{\varnothing}) - \min_x\{F(x)\}) + \frac{\mathcal{V}}{\sigma B}\Big) \enspace . \tag{6.2}$$

Recall we have chosen $B = \Theta(\frac{\mathcal{V}}{\varepsilon^2}) \geq 1$ so (6.2) implies

$$\mathbb{E}[\sigma^2\|\widehat{\mathsf{y}} - \widehat{\mathsf{y}}^*\|^2] \leq O\Big(\frac{\sigma^{1/3}L^{2/3}B^{2/3}\Delta_F}{T} + \varepsilon^2\Big) \leq O\Big(\frac{\sigma^{1/3}L^{2/3}\mathcal{V}^{2/3}\Delta_F}{\varepsilon^{4/3}T} + \varepsilon^2\Big)$$

In other words, as long as $T \geq \Omega\Big(\frac{(L^2\sigma)^{1/3}\Delta_F\cdot\mathcal{V}^{2/3}}{\varepsilon^{10/3}}\Big)$, we have $\mathbb{E}[\sigma^2\|\widehat{\mathsf{y}} - \widehat{\mathsf{y}}^*\|^2] \leq O(\varepsilon^2)$.

If we use $\mathtt{SGD3}^{\mathsf{sc}}$ of Theorem 5.7 to minimize the convex function $G(x) \stackrel{\text{def}}{=} F(x) + \sigma\|x - \widehat{\mathsf{y}}\|^2$ starting from $x = \widehat{\mathsf{x}}$, we get an output $x^{\mathsf{out}}$ satisfying[18]

$$\begin{aligned}\mathbb{E}[\|\mathcal{G}_{F^s,\eta}(x^{\mathsf{out}})\|] &\leq O\big(\frac{\sqrt{\mathcal{V}}}{\sqrt{T_{\mathsf{sgd}}}}\log^{3/2}\frac{L}{\sigma}\big) + \big(1 - \frac{\sigma}{L}\big)^{\Omega(T_{\mathsf{sgd}}/\log(L/\sigma))}\sigma\mathbb{E}[\|\widehat{\mathsf{y}} - \widehat{\mathsf{y}}^*\|] \\ &\leq O\big(\frac{\sqrt{\mathcal{V}}}{\sqrt{T_{\mathsf{sgd}}}}\log^{3/2}\frac{L}{\sigma}\big) + \big(1 - \frac{\sigma}{L}\big)^{\Omega(T_{\mathsf{sgd}}/\log(L/\sigma))}\varepsilon \enspace .\end{aligned}$$

[18]More specifically, we apply Theorem 5.7 for $G(x) = \psi(x) + \frac{1}{n}\sum_{i\in[n]}\big(f_i(x) + \sigma\|x - \widehat{\mathsf{x}}\|^2\big)$. It satisfies Problem (5.2) with the same smoothness $O(L)$, the same strong convexity $\sigma$, and the same variance bound $\mathcal{V}$.

In other words, as long as $T_{\mathsf{sgd}} \geq \Omega\Big(\frac{\mathcal{V}}{\varepsilon^2}\log^3\frac{L}{\sigma}\Big) \geq \Omega\Big(\frac{L}{\sigma}\log\frac{L}{\sigma}\Big)$, we have $\mathbb{E}[\|\mathcal{G}_{F^s,\eta}(x^{\mathsf{out}})\|] \leq O(\varepsilon)$. Finally, since $F^s(x) = F(x) + \sigma\|x - \widehat{\mathsf{x}}^s\|^2$ satisfies the assumption of $G(x)$ in Lemma 5.5, applying Lemma 5.5, we conclude that $\mathbb{E}[\|\mathcal{G}_{F,\eta}(x^{\mathsf{out}})\|] \leq O(1)\cdot\big(\mathbb{E}[\|\mathcal{G}_{F^s,\eta}(x^{\mathsf{out}})\|]+\mathbb{E}[\sigma\|\widehat{\mathsf{y}}-\widehat{\mathsf{y}}^*\|]\big) \leq O(\varepsilon)$. □

# 7 Natasha 2: Finding Local Minima

In this section, we study the problem finding approximate local minimum for

$$f(x) \stackrel{\text{def}}{=} \tfrac{1}{n}\textstyle\sum_{i=1}^n f_i(x) \ , \tag{7.1}$$

where

1. each $f_i(x)$ is possibly nonconvex but $L$-smooth,
2. the average $f(x)$ is possibly nonconvex, but second-order smooth with parameter $L_2$, and
3. the stochastic gradients $\nabla f_i(x)$ have a bounded variance, that is

$$\forall x \in \mathbb{R}^d\colon \quad \mathbb{E}_{i\in_R[n]}\|\nabla f(x) - \nabla f_i(x)\|^2 \leq \mathcal{V} \ .$$

This is the exact same setting studied by offline methods [1, 13, 14] and by online method SGD [22], except that the results in [1, 13, 14] did not assume any bound on variance. (Recall that variance bound is only necessary for online methods, see Table 1.) [19]

**Algorithm.** Our pseudocode $\texttt{Natasha2}^{\mathsf{full}}$ is given in Algorithm 5. It starts from a vector $y_0 \in \mathbb{R}^d$ and is divided into iterations $k = 0, 1, \dots$. In each iteration $k$, it *either* finds a vector $v \in \mathbb{R}^d$ such that $v^\top \nabla^2 f(y_k) v \leq -\frac{\delta}{2}$, *or* conclude that $\nabla^2 f(y_k) \succeq -\delta I$. This can be done via Oja's algorithm in Section 5.1.

- If $v^\top \nabla^2 f(y_k) v \leq -\frac{\delta}{2}$, we choose $y_{k+1} \leftarrow y_k + \frac{\delta}{L_2}v$ and $y_{k+1} \leftarrow y_k - \frac{\delta}{L_2}v$ each with probability 1/2. We call this a second-order step.
- If $\nabla^2 f(y_k) \succeq -\delta I$, then we define $F(x) = F^k(x) \stackrel{\text{def}}{=} f(x) + L(\max\{0, \|x - y_k\| - \frac{\delta}{L_2}\})^2$, and apply $\texttt{Natasha1.5}^{\mathsf{full}}$ for one full epoch (i.e., $T' = 1$). We call this a first-order step.

  Note that $\texttt{Natasha1.5}^{\mathsf{full}}$ returns two points $\widehat{\mathsf{y}}$ and $y^+$. We move to $y_{k+1} \leftarrow y^+$.

Finally, we terminate $\texttt{Natasha2}^{\mathsf{full}}$ whenever $N_1$ iterations of *first-order* steps are met. We select a random $\widehat{\mathsf{y}}$ along the $N_1$ first-order steps, and prune it using convex SGD. This is similar to the pruning step of $\texttt{Natasha1.5}^{\mathsf{full}}$.

Recall that $F(x)$ is $5L$-smooth and $3\delta$-strongly convex (see Claim 7.2). Thus, when applying $\texttt{Natasha1.5}^{\mathsf{full}}$, we can choose smoothness parameter $\widetilde{L}$ and bounded nonconvexity parameter $\widetilde{\sigma}$ for any $\widetilde{L} \geq 5L$ and $\widetilde{\sigma} \geq 3\delta$. Unfortunately, technical difficulties prevent us from always choosing $\widetilde{L} = 5L$ and $\widetilde{\sigma} = 3\delta$. We specify in Line 4 of $\texttt{Natasha2}^{\mathsf{full}}$ some special values for $\widetilde{L}$ and $\widetilde{\sigma}$ in order to tackle some boundary cases.

*Remark* 7.1. For instance, to provide a good control on the distance $\|x - y_k\|$ (see Section 3.2), we sometimes have to increase $\widetilde{\sigma}$ so that the distance $\|x - y_k\|$ becomes smaller (recall $\texttt{Natasha1.5}^{\mathsf{full}}$ performs retraction with weight $\widetilde{\sigma}$; so the larger $\widetilde{\sigma}$ is, the smaller $\|x - y_k\|$ becomes).

[19] Like in [1, 13, 14, 22], we do not include the proximal term $\psi(\cdot)$ when finding approximate local minima, because it can be tricky to define what local minima mean when $\psi(\cdot)$ is present.

**Algorithm 5** $\texttt{Natasha2}^{\mathsf{full}}(f, y_0, \varepsilon, \delta)$

**Input:** function $f(x)$ satisfying Problem (7.1), starting vector $y_0$, target accuracy $\varepsilon > 0$ and $\delta > 0$.
◇ *assume* $\mathcal{V} \geq \Omega(\varepsilon^2)$ *and* $\delta^4 \leq O(\mathcal{V}\varepsilon L_2^3/L^2)$

**if** $L_2 \geq \frac{L\delta}{\mathcal{V}^{1/3}\varepsilon^{1/3}}$ **then** ◇ *the boundary case for large* $L_2$
    $\widetilde{L} = \widetilde{\sigma} \leftarrow \Theta(\frac{L_2 \mathcal{V}^{1/3}\varepsilon^{1/3}}{\delta}) \in [L, \infty)$.
**else** ◇ *the most interesting case*
    $\widetilde{L} \leftarrow L$ and $\widetilde{\sigma} \leftarrow \Theta\big(\max\{\frac{\mathcal{V}\varepsilon L_2^3}{L^2\delta^3}, \frac{\varepsilon L}{\mathcal{V}^{1/2}}\}\big) \in [\delta, L]$.
**end if**
$B \leftarrow \Theta(\mathcal{V}/\varepsilon^2)$; $p \leftarrow \Theta\big((\frac{\widetilde{\sigma}^2}{\widetilde{L}^2}B)^{1/3}\big)$ $\alpha \leftarrow \Theta(\frac{\widetilde{\sigma}}{p^2\widetilde{L}^2})$. ◇ *same* $p$ *and* $\alpha$ *as* $\texttt{Natasha1.5}^{\mathsf{prune}}$
$X \leftarrow []$; $N_1 \leftarrow \Theta(\frac{\widetilde{\sigma}\Delta_f}{p\varepsilon^2})$, where $\Delta_f$ is an upper bound on $f(y_0) - \min_y\{f(y)\}$.
**for** $k \leftarrow 0$ **to** $\infty$ **do**
    Apply Oja's algorithm to find minEV of $\nabla^2 f(y_k)$. ◇ *use Lemma 5.3 with* $T_{\mathsf{oja}} = \Theta(\frac{L^2}{\delta^2}\log(dk))$

2nd-order step:
    **if** $v \in \mathbb{R}^d$ is found s.t. $v^\top \nabla^2 f(y_k) v \leq -\frac{\delta}{2}$ **then**
        $y_{k+1} \leftarrow y_k \pm \frac{\delta}{L_2}v$ where the sign is random.
    **else** ◇ *it satisfies* $\nabla^2 f(y_k) \succeq -\delta\mathbf{I}$ *w.p.* $\geq 1 - \frac{1}{20(k+1)^2}$*, see Lemma 5.3.*

1st-order step:
        $F(x) = F^k(x) \stackrel{\text{def}}{=} f(x) + L(\max\{0, \|x - y_k\| - \frac{\delta}{L_2}\})^2$.
        ◇ $F(\cdot)$ *is* $O(\widetilde{L})$*-smooth and* $O(\widetilde{\sigma})$*-nonconvex*
        $(\widehat{y}_k, y_{k+1}) \leftarrow \texttt{Natasha1.5}^{\mathsf{full}}(F, y_k, B, p, 1, \alpha)$
        $X \leftarrow [X, (y_k, \widehat{y}_k)]$.
        break the for loop if have performed $N_1$ first-order steps.
    **end if**
**end for**
$(y, \widehat{y}) \leftarrow$ a random pair in $X$. ◇ *in practice, letting* $x^{\mathsf{out}} =$ *the last* $\widehat{y}_k$ *should be good enough*
define $G(x) \stackrel{\text{def}}{=} f(x) + L(\max\{0, \|x - y\| - \frac{\delta}{L_2}\})^2 + \widetilde{\sigma}\|x - \widehat{y}\|^2$. ◇ $G(x)$ *is* $\widetilde{\sigma}$*-strongly convex*
$x^{\mathsf{out}} \leftarrow \texttt{SGD3}^{\mathsf{sc}}(G, \widehat{y}, \widetilde{\sigma}, O(\widetilde{L}), T_{\mathsf{sgd}})$. ◇ *with* $T_{\mathsf{sgd}} = \Theta(\frac{\mathcal{V}}{\varepsilon^2}\log^3\frac{\widetilde{L}}{\widetilde{\sigma}})$
**return** $x^{\mathsf{out}}$.

## 7.1 Natasha 2: Auxiliary Claims

**Claim 7.2.** *If $f(x)$ is $L$-smooth and second-order smooth with parameter $L_2$, and $y \in \mathbb{R}^d$ is a point such that $\nabla^2 f(y) \succeq -\delta\mathbf{I}$ for some $\delta > 0$, then the function*

$$F(x) = f(x) + L\big(\max\{0, \|x - y\| - \tfrac{\delta}{L_2}\}\big)^2$$

*is $5L$ smooth and $3\delta$-nonconvex.*

*Proof.* This is a simple consequence of the smoothness definition, see proofs in [13, Lemma 4.1]. □

**Claim 7.3.** *If $v^\top \nabla^2 f(y_k) v \leq -\frac{\delta}{2}$ and we run a second-order step, then $f(y_k) - \mathbb{E}[f(y_{k+1})] \geq \frac{\delta^3}{12L_2^2}$.*

*Proof.* Suppose $y_{k+1} = y_k \pm \eta v$ where $\|v\| = 1$ and $\eta = \frac{\delta}{L_2}$, then by the second-order smoothness,

$$\begin{aligned} f(y_k) - \mathbb{E}[f(y_{k+1})] &\geq \mathbb{E}\big[\langle \nabla f(y_k), y_k - y_{k+1}\rangle - \frac{1}{2}(y_k - y_{k+1})^\top \nabla^2 f(y_k)(y_k - y_{k+1}) - \frac{L_2}{6}\|y_k - y_{k+1}\|^3\big] \\ &= -\frac{\eta^2}{2} v^\top \nabla^2 f(y_k) v - \frac{L_2\eta^3}{6}\|v\|^3 \geq \frac{\eta^2\delta}{4} - \frac{L_2\eta^3}{6} = \frac{\delta^3}{12L_2^2} \ . \end{aligned}$$

□

**Claim 7.4.** *If $\nabla^2 f(y_k) \succeq -\delta \mathbf{I}$ and we run a first-order step, then by Lemma 6.5,*

$$f(y_k) - \mathbb{E}[f(y_{k+1})] \geq \frac{1}{4}\mathbb{E}\Big[\frac{\widetilde{\sigma}}{p}\|y_k - \widehat{y}_k\|^2 + \widetilde{\sigma}p\|\widehat{y}_k - \widehat{y}_k^*\|^2\Big] - \frac{3p\mathcal{V}}{2\widetilde{\sigma}B} \enspace,$$

*where* $\widehat{\mathsf{y}}_k^* \stackrel{\mathrm{def}}{=} \arg\min_x\{F^k(x) + \widetilde{\sigma}\|x - \widehat{\mathsf{y}}_k\|^2\}$.

*Proof.* We can apply Lemma 6.5 with $L = O(\widetilde{L})$ and $\sigma = O(\widetilde{\sigma})$, because $F^k(x)$ is $O(\widetilde{L})$-smooth and $O(\widetilde{\sigma})$-nonconvex (see Claim 7.2) and recall $\widetilde{\sigma} \geq \delta$. Also, we always choose $p \stackrel{\mathrm{def}}{=} \big(\frac{\widetilde{\sigma}^2}{48\widetilde{L}^2}B\big)^{1/3}$, $m = B/p$, and $\alpha \stackrel{\mathrm{def}}{=} \frac{8}{\widetilde{\sigma}m} = \frac{\widetilde{\sigma}}{6p^2\widetilde{L}^2} \leq \frac{\widetilde{\sigma}}{6\widetilde{L}^2} \leq \frac{1}{2\widetilde{L}+4\widetilde{\sigma}}$. These parameters satisfy the prerequisite of Lemma 6.5. Since $\widehat{\mathsf{y}} = \widehat{\mathsf{x}}^s$ where $s$ is a random subepoch $s \in \{0, 1, \dots, p-1\}$ in $\texttt{Natasha1.5}^{\mathsf{full}}$, we have

$$\begin{aligned}
&\mathbb{E}\Big[\frac{\widetilde{\sigma}}{2p}\|y_k - \widehat{y}_k\|^2 + \frac{\widetilde{\sigma}p}{4}\|\widehat{y}_k - \widehat{y}_k^*\|^2\Big] = \mathbb{E}\Big[\frac{\widetilde{\sigma}}{2p}\|\widehat{\mathsf{x}}^0 - \widehat{\mathsf{x}}^s\|^2 + \frac{\widetilde{\sigma}p}{4}\|\widehat{\mathsf{x}}^s - x_*^s\|^2\Big] \\
&\leq \mathbb{E}\Big[\frac{\widetilde{\sigma}}{2}\sum_{s=0}^{p-1}\Big(\|\widehat{\mathsf{x}}^s - \widehat{\mathsf{x}}^{s+1}\|^2 + \frac{1}{2}\|\widehat{\mathsf{x}}^s - x_*^s\|^2\Big)\Big] \overset{①}{\leq} \mathbb{E}\Big[F(\widehat{\mathsf{x}}^0) - F(\widehat{\mathsf{x}}^p)\Big] + \frac{3p\mathcal{V}}{2\widetilde{\sigma}B} \enspace.
\end{aligned}$$

Above, inequality ① uses Lemma 6.5. Finally, we note that $F(\widehat{\mathsf{x}}^0) = F(y_k) = f(y_k)$ and $F(\widehat{\mathsf{x}}^p) = f(\widehat{\mathsf{x}}^p) + L\big(\max\{0, \|\widehat{\mathsf{x}}^p - y\| - \frac{\delta}{L_2}\}\big)^2 \geq f(\widehat{\mathsf{x}}^p) = f(y_{k+1})$, so finish the proof. □

## 7.2 Natasha 2: Main Theorem

We state the main theorem of $\texttt{Natasha2}^{\mathsf{full}}$ as follows.

**Theorem 2.** *Consider Problem (7.1) with a starting vector $y_0$. For any $\varepsilon > 0$ and $\delta \in (0, L]$, under assumptions $\mathcal{V} \geq \Omega(\varepsilon^2)$ and $\delta^4 \leq O(\mathcal{V}\varepsilon L_2^3/L^2)$, the output $x^{\mathsf{out}} = \texttt{Natasha2}^{\mathsf{full}}(f, y_0, \varepsilon, \delta)$ satisfies, with probability at least 2/3,*

$$\|\nabla f(x^{\mathsf{out}})\| \leq \varepsilon \quad \text{and} \quad \nabla^2 f(x^{\mathsf{out}}) \succeq -3\delta\mathbf{I} \enspace.$$

*The total gradient complexity $T$ is*

$$T = \widetilde{O}\Big(\frac{\mathcal{V}^2}{\varepsilon^2} + \frac{L_2^2 L^2 \Delta_f}{\delta^5} + \frac{L_2\Delta_f}{\varepsilon\delta}\cdot\Big(\frac{L^2}{\delta^2} + \frac{\mathcal{V}}{\varepsilon^2}\Big) + \frac{L\Delta_f}{\varepsilon\mathcal{V}^{1/2}}\cdot\frac{L^2}{\delta^2}\Big)$$

*Above, $\Delta_f$ is any known upper bound on $f(y_0) - \min_y\{f(y)\}$.*

*Remark* 7.5. In practice, one can just choose $N_1$, the number of first-order updates in $\texttt{Natasha2}^{\mathsf{full}}$, as sufficiently large, without the necessity of knowing $\Delta_f$.

*Remark* 7.6. As a sanity check, our formula for $T$ in Theorem 2 is scaling invariant: if $f(x)$ increases by a factor $C$, then $\Delta_f$, $L$, $\varepsilon$, and $L_2$ each increases by $C$, and $\mathcal{V}$ increases by $C^2$.

**Corollary 7.7.** *If we assume $L, L_2, \Delta_f$ and $\mathcal{V}$ are constants, then $\texttt{Natasha2}^{\mathsf{full}}$ finds $x^{\mathsf{out}}$ satisfying*

$$\|\nabla f(x^{\mathsf{out}})\| \leq \varepsilon \quad \text{and} \quad \nabla^2 f(x^{\mathsf{out}}) \succeq -\delta\mathbf{I}$$

*in gradient complexity $T = \widetilde{O}\big(\frac{1}{\delta^5} + \frac{1}{\delta\varepsilon^3}\big)$ for every $\delta \in (0, \varepsilon^{1/4}]$. Since when $\delta > \varepsilon^{1/4}$, we can replace $\delta$ with $\varepsilon^{1/4}$, this complexity can be re-written as $T = \widetilde{O}\big(\frac{1}{\delta^5} + \frac{1}{\varepsilon^{3.25}} + \frac{1}{\delta\varepsilon^3}\big)$.*

## 7.3 Natasha 2: Proof of Theorem 2

We use the big-$\Theta$ notion to hide absolute constants, in order to simplify notations.

*Proof of Theorem 2.* Recall $N_1 = \Theta\big(\frac{\widetilde{\sigma}\Delta_f}{p\varepsilon^2}\big)$ is the number of first-order steps. We denote by $N_2$ the actual number of second-order steps, which is a *random variable*.

We first note that each call of Oja's algorithm succeeds with probability at least $1 - \frac{1}{20(k+1)^2}$, and therefore by $\sum_{k=1}^{\infty} k^{-2} < 1.65$, with probability at least $1 - \frac{1}{12}$ (over the randomness of Oja's algorithm only), all occurrences of Oja's algorithm succeed. In the remainder of the proof, we shall always assume that this event happens. In other words, in Line 9 of $\mathtt{Natasha2}^{\mathsf{full}}$, it either finds $v^\top \nabla^2 f(y_k) v \leq -\frac{\delta}{2}$ or if not, conclude that $\nabla^2 f(y_k) \succeq -\delta \mathbf{I}$. (Recall Lemma 5.3.)

Let us define random variables $\Delta_1, \Delta_2$ the total amount of objective decrease during first-order and second-order steps respectively.[20] By Claim 7.4 and the fact that there are exactly $N_1$ first-order steps, we have $\mathbb{E}[\Delta_1] \geq -\frac{3p\mathcal{V}}{2\widetilde{\sigma}B} N_1 = -\Theta(\frac{\mathcal{V}}{B\varepsilon^2}) \cdot \Delta_f \geq -\Delta_f$, where we have used $N_1 = \Theta\big(\frac{\widetilde{\sigma}\Delta_f}{p\varepsilon^2}\big)$ and $B = \Theta(\frac{\mathcal{V}}{\varepsilon^2})$.

**Accuracy.** Since $\Delta_1 + \Delta_2 \leq \Delta_f$ and $\mathbb{E}[\Delta_2] \geq 0$ by Claim 7.3, we conclude that if we select $k = 0, 1, \dots,$ at random among the $N_1$ first-order steps, then

$$\mathbb{E}[f(y_k) - f(y_{k+1})] \leq \frac{\mathbb{E}[\Delta_1]}{N_1} \leq \frac{\Delta_f - \mathbb{E}[\Delta_2]}{N_1} \leq \frac{\Delta_f}{N_1} \enspace .$$

Denote by $y = y_k$, $\widehat{\mathsf{y}} = \widehat{\mathsf{y}}_k$, and $\widehat{\mathsf{y}}^* = \arg\min_x \{F^k(x) + \widetilde{\sigma}\|x - \widehat{\mathsf{y}}_k\|^2\}$ for this random choice of $k$. Combining $\mathbb{E}[f(y_k) - f(y_{k+1})] \leq \frac{\Delta_f}{N_1}$ and Claim 7.4, we have

$$\mathbb{E}\Big[\frac{\widetilde{\sigma}}{p}\|y - \widehat{y}\|^2 + \widetilde{\sigma}p\|\widehat{y} - \widehat{y}^*\|^2\Big] \leq O\Big(\frac{\Delta_f}{N_1} + \frac{p\mathcal{V}}{\widetilde{\sigma}B}\Big) = O\Big(\frac{p\varepsilon^2}{\widetilde{\sigma}}\Big) \enspace .$$

By Markov's bound, with probability at least, $1 - \frac{1}{12}$, we have

$$\frac{\widetilde{\sigma}}{p}\|y - \widehat{y}\|^2 + \widetilde{\sigma}p\|\widehat{y} - \widehat{y}^*\|^2 \leq O\Big(\frac{p\varepsilon^2}{\widetilde{\sigma}}\Big) \enspace . \tag{7.2}$$

Now, recall that

$$F(x) \stackrel{\mathrm{def}}{=} f(x) + L(\max\{0, \|x - y\| - \tfrac{\delta}{L_2}\})^2 \quad \text{and} \quad G(x) = F(x) + \widetilde{\sigma}\|x - \widehat{\mathsf{y}}\|^2$$

we can apply $\mathtt{SGD3}^{\mathsf{sc}}$ for gradient complexity $T_{\mathsf{sgd}}$ to minimize $G(x)$. Let the output be $x^{\mathsf{out}}$. Using Theorem 5.7, we have with probability at least $1 - \frac{1}{12}$[21]

$$\|\nabla G(x^{\mathsf{out}})\|^2 \leq O\Big(\frac{\mathcal{V}}{T_{\mathsf{sgd}}}\log^3\frac{\widetilde{L}}{\widetilde{\sigma}}\Big) + \big(1 - \frac{\widetilde{\sigma}}{\widetilde{L}}\big)^{\Omega(T_{\mathsf{sgd}}/\log(\widetilde{L}/\widetilde{\sigma}))}\widetilde{\sigma}^2\|\widehat{\mathsf{y}} - \widehat{\mathsf{y}}^*\|^2 \enspace . \tag{7.3}$$

Using Lemma 5.5, we have

$$\|\nabla F(x^{\mathsf{out}})\|^2 + \widetilde{\sigma}^2\|x^{\mathsf{out}} - \widehat{\mathsf{y}}\|^2 \leq O\big(\widetilde{\sigma}^2\|\widehat{\mathsf{y}}^* - \widehat{\mathsf{y}}\|^2 + \|\nabla G(x^{\mathsf{out}})\|^2\big) \enspace . \tag{7.4}$$

Combining (7.2), (7.3), and (7.4), and our choice $T_{\mathsf{sgd}} = \Theta\big(\frac{\mathcal{V}}{\varepsilon^2}\log^3\frac{\widetilde{L}}{\widetilde{\sigma}}\big) \geq \Theta\big(\frac{\widetilde{L}}{\widetilde{\sigma}}\log\frac{\widetilde{L}}{\widetilde{\sigma}}\big)$ we have

$$\|\nabla F(x^{\mathsf{out}})\|^2 + \widetilde{\sigma}^2\|x^{\mathsf{out}} - \widehat{\mathsf{y}}\|^2 \leq O\Big(\frac{p\varepsilon^2}{\widetilde{\sigma}} \cdot \frac{\widetilde{\sigma}^2}{\widetilde{\sigma}p} + \frac{\mathcal{V}}{T_{\mathsf{sgd}}}\log^3\frac{\widetilde{L}}{\widetilde{\sigma}}\Big) \leq \varepsilon^2 \enspace . \tag{7.5}$$

Recall that (in both cases, $L_2$ too large or not) we have chosen $\widetilde{\sigma}$ and $p$ so that

$$\frac{p\varepsilon}{\widetilde{\sigma}} \leq O\Big(\frac{\delta}{L_2}\Big) \enspace , \tag{7.6}$$

[20] More precisely, $\Delta_1 \stackrel{\mathrm{def}}{=} \sum_{k=0}^{\infty} \mathbb{I}[\text{iter } k \text{ exists and is a first-order step}] \cdot (f(y_k) - f(y_{k+1}))$ , and similarly for $\Delta_2$.

[21] More specifically, we apply Theorem 5.7 for $G(x) = \frac{1}{n}\sum_{i\in[n]}\big(f_i(x) + L(\max\{0, \|x - y\| - \frac{\delta}{L_2}\})^2 + \widetilde{\sigma}\|x - \widehat{\mathsf{y}}\|^2\big)$. It satisfies Problem (5.2) with the same smoothness $O(\widetilde{L})$, the same strong convexity $O(\widetilde{\sigma})$, and the same variance bound $\mathcal{V}$.

therefore, (7.2) implies

$$\|y-\widehat{y}\|^2 \leq O\big(\frac{p^2}{\widetilde{\sigma}^2}\varepsilon^2\big) \leq \big(\frac{\delta}{2L_2}\big)^2 \ . \tag{7.7}$$

By triangle inequality,

$$\|x^{\mathsf{out}}-y\| \leq \|x^{\mathsf{out}}-\widehat{\mathsf{y}}\| + \|\widehat{\mathsf{y}}-y\| \leq \frac{\varepsilon}{\widetilde{\sigma}} + \frac{\delta}{2L_2} \leq \frac{p\varepsilon}{\widetilde{\sigma}} + \frac{\delta}{2L_2} \leq \frac{\delta}{L_2} \ .$$

In other words, $x^{\mathsf{out}}$ is not too far away from $y$ and therefore by definition $F(x) \stackrel{\mathrm{def}}{=} f(x) + L(\max\{0, \|x-y\| - \frac{\delta}{L_2}\})^2$,

$$\nabla^2 F(x^{\mathsf{out}}) = \nabla^2 f(x^{\mathsf{out}}) \quad \text{and} \quad \nabla F(x^{\mathsf{out}}) = \nabla f(x^{\mathsf{out}}) \ .$$

This means $\nabla^2 f(x^{\mathsf{out}}) = \nabla^2 F(x^{\mathsf{out}}) \succeq -3\delta\mathbf{I}$ (by the $3\delta$-nonconvexity of $F(\cdot)$, see Claim 7.2) and $\|\nabla f(x^{\mathsf{out}})\| = \|\nabla F(x^{\mathsf{out}})\| \leq \varepsilon$ by (7.5). This finishes the proof of the accuracy of $\texttt{Natasha2}^{\mathsf{full}}$.

**Running Time.** Recall that random variable $N_2$ is the number of second-order steps. By Claim 7.3, we have

$$\mathbb{E}[N_2]\cdot\frac{\delta^3}{12L_2^2} \leq \mathbb{E}[\Delta_2] \leq \Delta_f - \mathbb{E}[\Delta_1] \leq 2\Delta_f \implies \mathbb{E}[N_2] \leq O\big(\frac{L_2^2\Delta_f}{\delta^3}\big) \ .$$

Therefore, with probability at least $1-\frac{11}{12}$, we have $N_2 \leq O\big(\frac{L_2^2\Delta_f}{\delta^3}\big)$. The remainder of the derivation always assumes this event happens.

The total gradient complexity $T$ consists of three parts:

- The gradient complexity for Oja's algorithms is at most $O\big((N_1+N_2)\frac{L^2}{\delta^2}\big)$.
- The gradient complexity for applying $\texttt{Natasha1.5}^{\mathsf{full}}$ for $N_1 = \Theta\big(\frac{\widetilde{\sigma}\Delta_f}{p\varepsilon^2}\big)$ times is at most $N_1\cdot B$.
- The gradient complexity for applying SGD in the end is $T_{\mathsf{sgd}} = O\big(\frac{\mathcal{V}}{\varepsilon^2}\log^3\frac{\widetilde{L}}{\widetilde{\sigma}}\big)$.

**Case 1.** Suppose $L_2 \geq \frac{L\delta}{\mathcal{V}^{1/3}\varepsilon^{1/3}}$. This corresponds to the case when $L_2$ is too large. Recall we have chosen $\widetilde{L} = \widetilde{\sigma} = \Theta(\frac{p\varepsilon L_2}{\delta}) = \Theta(\frac{L_2\mathcal{V}^{1/3}\varepsilon^{1/3}}{\delta}) \geq L$ and (7.6) is satisfied. The total gradient complexity is

$$\begin{aligned}T = \widetilde{O}\Big(T_{\mathsf{sgd}} + (N_1+N_2)\cdot\frac{L^2}{\delta^2} + N_1\cdot\frac{\mathcal{V}}{\varepsilon^2}\big) &\leq \widetilde{O}\Big(\frac{\mathcal{V}}{\varepsilon^2} + \big(\frac{\widetilde{L}\Delta_f}{\varepsilon^2 p} + \frac{L_2^2\Delta_f}{\delta^3}\big)\cdot\frac{L^2}{\delta^2} + \frac{\widetilde{L}\Delta_f}{\varepsilon^2 p}\cdot\frac{\mathcal{V}}{\varepsilon^2}\Big)\\ &\leq \widetilde{O}\Big(\frac{\mathcal{V}}{\varepsilon^2} + \big(\frac{L_2\Delta_f}{\varepsilon\delta} + \frac{L_2^2\Delta_f}{\delta^3}\big)\cdot\frac{L^2}{\delta^2} + \frac{L_2\Delta_f}{\varepsilon\delta}\cdot\frac{\mathcal{V}}{\varepsilon^2}\Big) \ .\end{aligned}$$

**Case 2.** Suppose $L_2 \leq \frac{L\delta}{\mathcal{V}^{1/3}\varepsilon^{1/3}}$. This is the *interesting case* and recall we have chosen $\widetilde{L} = L$ and $\widetilde{\sigma} = \Omega\big(\max\big\{\frac{\mathcal{V}\varepsilon L_2^3}{L^2\delta^3}, \frac{\varepsilon L}{\mathcal{V}^{1/2}}\big\}\big)$. One can verify that $\widetilde{\sigma} \in [\delta, L]$ and (7.6) is satisfied.

The total gradient complexity

$$\begin{aligned}T &= \widetilde{O}\Big(T_{\mathsf{sgd}} + (N_1+N_2)\cdot\frac{L^2}{\delta^2} + N_1\cdot\frac{\mathcal{V}}{\varepsilon^2}\big)\Big) \leq \widetilde{O}\Big(\frac{\mathcal{V}^2}{\varepsilon^2} + \big(\frac{\widetilde{\sigma}\Delta_f}{\varepsilon^2 p} + \frac{L_2^2\Delta_f}{\delta^3}\big)\cdot\frac{L^2}{\delta^2} + \frac{\widetilde{\sigma}\Delta_f}{\varepsilon^2 p}\cdot\frac{\mathcal{V}}{\varepsilon^2}\Big)\\ &= \widetilde{O}\Big(\frac{\mathcal{V}^2}{\varepsilon^2} + \frac{L_2^2L^2\Delta_f}{\delta^5} + \Big(\frac{L^{2/3}\widetilde{\sigma}^{1/3}\Delta_f}{\varepsilon^{4/3}\mathcal{V}^{1/3}}\Big)\cdot\Big(\frac{L^2}{\delta^2} + \frac{\mathcal{V}}{\varepsilon^2}\Big)\Big)\\ &= \widetilde{O}\Big(\frac{\mathcal{V}^2}{\varepsilon^2} + \frac{L_2^2L^2\Delta_f}{\delta^5} + \frac{L_2\Delta_f}{\varepsilon\delta}\cdot\Big(\frac{L^2}{\delta^2} + \frac{\mathcal{V}}{\varepsilon^2}\Big) + \frac{L\Delta_f}{\varepsilon\mathcal{V}^{1/2}}\cdot\frac{L^2}{\delta^2}\Big) \ . \qquad \square\end{aligned}$$

## Acknowledgements

We would like to thank Lin Xiao for suggesting reference [47, Lemma 3.7], and Yurii Nesterov for useful discussions on the convex version of this problem, Sébastien Bubeck, Yuval Peres, and Lin Xiao for discussing notations, and Chi Jin for discussing reference [46].

# Appendix

| | algorithm | gradient complexity $T$ | variance bound | Lipschitz smooth | 2nd-order smooth |
|---|---|---|---|---|---|
| convex only | `GD` [35] | $O(n\varepsilon^{-1})$ ♭ | no | needed | no |
| | `AccGD` [35] | $\widetilde{O}(n\varepsilon^{-1/2})$ ♭ | no | needed | no |
| | `SVRG` [28, 50] | $\widetilde{O}(n+\varepsilon^{-1})$ ♭ | no | needed | no |
| | `AccSVRG` [20, 32] or `Katyusha` [2] | $\widetilde{O}(n+n^{1/2}\varepsilon^{-1/2})$ | no | needed | no |
| approximate stationary points | `GD` (folklore [33]) | $O(n\varepsilon^{-2})$ ♭ (see Appendix B) | no | needed | no |
| | `SVRG` [6, 39] | $O(n+n^{2/3}\varepsilon^{-2})$ ♭ | no | needed | no |
| | `repeatSVRG` [3, 13] | $O(n+n^{3/4}\sigma^{1/2}\varepsilon^{-2}+n\sigma\varepsilon^{-2})$ | no | needed | no |
| | `Natasha1` [3] | $O(n+n^{1/2}\varepsilon^{-2}+\sigma^{1/3}n^{2/3}\varepsilon^{-2})$ | no | needed | no |
| | CHDS [14] | $\widetilde{O}(n\varepsilon^{-1.75})$ | no | needed | needed |
| approximate local minima | perturbed `GD` [27] | $\widetilde{O}(n\varepsilon^{-2})$ ♭ | no | needed | needed |
| | CHDS [13] `FastCubic` [1] | $\widetilde{O}(n\varepsilon^{-1.75})$ ♭ | no | needed | needed |
| | [10, 40] | $\widetilde{O}(n+n^{2/3}\varepsilon^{-2})$ | no | needed | needed |
| | CHDS [13] `FastCubic` [1] | $\widetilde{O}(n\varepsilon^{-1.5}+n^{3/4}\varepsilon^{-1.75})$ | no | needed | needed |

Table 2: Comparison of ***offline*** methods for finding $\|\nabla f(x)\| \leq \varepsilon$. This table is for ***reference purpose only***. Following tradition, in these complexity bounds, we assume variance and smoothness parameters as constants, and only show the dependency on $n, d, \varepsilon$ and the bounded nonconvexity parameter $\sigma \in (0,1)$. We use ♭ to indicate the result is outperformed. Note that $n+n^{1/2}\varepsilon^{-1/2} \leq O(n+\varepsilon^{-1})$ so `SVRG` is outperformed by `AccSVRG`/`Katyusha` in the convex case.

## A Other Related Works

**Variance Reduction.** Methods based on variance reduction were first introduced for convex optimization. The first such method is SAG [41], but SAG cannot handle proximal terms so cannot be applied to tasks such as Lasso, SVM, etc. This was later fixed in two ways: the `SVRG` approach we adopted in this paper [28, 50], and the SAGA approach we did not use [18]. In the convex setting, variance reduction can be made accelerated using momentum [2].

The first "nonconvex use" of variance reduction is by Shalev-Shwartz [44], who assumes that each $f_i(x)$ is nonconvex but their average $f(x)$ is still convex. This result can also be made accelerated using momentum [4].

The first truly nonconvex use of variance reduction (i.e., for $f(x)$ being also nonconvex) is independently by Reddi et al. [39] and Allen-Zhu and Hazan [6], both in March 2016. Later, Lei et al. [30] made this approach online by their SCSG method.

All of these cited methods except SCSG are offline.

**Second-Order Methods.** If one is allowed to invert the Hessian matrix, then cubic-regularized Newton's method [34] converges in $1/\varepsilon^{3/2}$ iterations. Since its per-iteration complexity is very high, we have not included it in Table 2. Agarwal et al. [1] showed that the same cubic-regularized Newton's method can be implemented using only $T = \widetilde{O}\big(\frac{n}{\varepsilon^{1.5}} + \frac{n^{3/4}}{\varepsilon^{1.75}}\big)$ computations of stochastic gradients and Hessian-vector products. A similar result can also be obtained via the concurrent work of Carmon et al. [13]. These cited methods are offline, but are later turned into online ones by Tripuraneni et al. [46].

**Stochastic Eigenvector Computations.** The problem of finding the leading $k$ eigenvectors for a matrix $\mathbf{M} = \frac{1}{n}\sum_{i=1}^{n} \mathbf{M}_i$ (say, to an error $\delta > 0$) has received lots of attention in machine learning and theoretical computer science. In the *offline* setting, one can apply both variance reduction and acceleration techniques to achieve the fastest convergence rate $\delta^{-1/2}$. The first such result for $k = 1$ was [21] and for $k > 1$ was [7]. In the *online* setting, sampling lower bound prevents us from using variance reduction or acceleration, so the optimal convergence rate is $\delta^{-2}$ (see [9]). In this regime, Oja's algorithm can be viewed as a simple online stochastic version of power method, and achieves optimal complexity for both $k = 1$ [9], and for $k > 1$ at least when matrices $\mathbf{M}_i$ are rank-one [8].

**One-point Convexity.** Lots of recent progresses in nonconvex machine learning were based on showing that, if the data is sufficiently random, then the nonconvex function $f(x)$ satisfies for instance $\langle \nabla f(x), x - x^* \rangle \geq \Omega(\|x - x^*\|^2)$ or $\|\nabla f(x)\|^2 \geq \Omega(1) \cdot (f(x) - f(x^*))$. This is what we summarize as "one-point convexity" because it asks for a weak version of convexity between any point $x$ and the global minimum $x^*$ (where $x^*$ is assumed to exist). One-point convex functions are easy to minimize: for instance, gradient descent always converges to the global minimum.

However, one-point convexity is only known to apply to relatively simpler nonconvex tasks such as matrix completion [45], dictionary learning [12], phase retrieval [15] and a two-layer neural network [31], but not for complicated tasks such as training a deep neural network.

**Heuristics for Nonconvex Optimization.** Experimentalists have used AdaGrad [19], AdaDelta [49], Adam [29], and many other variants of SGD to train neural networks faster. For instance, AdaGrad applies a diagonal matrix to precondition (thus re-scale) the coordinates of $f(x)$. This is effective for neural networks, because weight variables $x_i$ across different layers of the network should be trained using separate step lengths. AdaDelta is built on AdaGrad but calculates the step length based on a window of accumulated past gradients. To the best of our knwoledge, there is no theoretical evidence that preconditioning methods like AdaGrad or AdaDelta affect the *convergence rate* of SGD in the nonconvex setting. Adam is similar to AdaDelta, but it adds Nesterov's momentum [33] on the top. To the best of our knowledge, there is no theoretical evidence that Nesterov's momentum helps improve the *convergence rate* of SGD for nonconvex functions (unless one imposes strong assumptions such as one-point convexity).

Neural network algorithms using Hessian-vector products have received some attention by experimentalists as well, see for instance [26] and the references therein. Such methods are referred to as Hessian-free methods. To the best of our knowledge, there is no theoretical evidence that they can improve the *convergence rate* of SGD for nonconvex functions.

### A.1 Other Extensions

**Mini-Batch.** Just like most stochastic methods, our $\texttt{Natasha1.5}^{\mathsf{full}}$ and $\texttt{Natasha2}^{\mathsf{full}}$ also have their mini-batch variants with provably convergence, which can be implemented via parallel computations and thus be applicable to even larger scales of machine learning tasks. In particular, whenever a gradient $\nabla f_i(x)$ is computed in $\texttt{Natasha1.5}^{\mathsf{full}}$, one can use $\frac{1}{|S|}\sum_{i\in S} f_i(x)$ instead for a random mini-batch $S \subseteq [n]$; whenever a Hessian-vector product $\nabla^2 f_i(x)\cdot v$ is needed in $\texttt{Natasha2}^{\mathsf{full}}$, one can replace it with $\frac{1}{|S|}\sum_{i\in S}\nabla^2 f_i(x)\cdot v$ for a random mini-batch $S$ as well. All of our theorems can be restated in such settings, but we refrain from doing so in order to keep the notations simple.

**Strict-Saddle Functions.** Some recent results [13, 22, 27] also state their convergence theorems using the strict-saddle language. These are just corollaries of finding $\varepsilon$-approximate local minima. For instance, in [27], a function $f(x)$ is $(\theta,\varepsilon,\delta)$-strict saddle if for any point $x\in\mathbb{R}^d$, one of the following three holds: (1) $\|\nabla f(x)\| > \varepsilon$, (2) $\lambda_{\min}(\nabla^2 f(x)) < -\delta$, or (3) $x$ is $\theta$ close to an *exact* local minimum. By applying Theorem 2, our $\texttt{Natasha2}^{\mathsf{full}}$ is able to find a point $\theta$ closes to an exact local minimum in gradient complexity $T = \widetilde{O}\big(\frac{1}{\delta^5} + \frac{1}{\delta\varepsilon^3} + \frac{1}{\varepsilon^{3.25}}\big)$.

# B Convergence of GD and SGD for Nonconvex Functions

We are not sure what is the earliest reference for showing that gradient descent and stochastic gradient descent converge to approximate stationary points. Both statements are simple to prove.

If a function $f(x)$ is $L$-smooth, then classical gradient descent theory (cf. textbook by Nesterov [33]) shows that, if we iteratively update $x_{t+1}\leftarrow x_t - \eta\nabla f(x_t)$ for step length $\eta = \frac{1}{L}$, then

$$f(x_t) - f(x_{t+1}) \geq \langle \nabla f(x_t), x_t - x_{t+1}\rangle - \frac{L}{2}\|x_t - x_{t+1}\|^2 = \frac{1}{2L}\|\nabla f(x_t)\|^2 \ .$$

Therefore, if we perform $N$ gradient updates $t = 0,1,\dots,N-1$, there must exist some point $x_t$ satisfying $\|\nabla f(x_t)\|^2 \leq O\big(\frac{L(f(x_0)-f(x^*))}{N}\big)$. Since each gradient computation $\nabla f(x)$ requires one to compute $n$ individual $\nabla f_i(x)$, this totals to gradient complexity $T = Nn \propto \frac{n}{\varepsilon^2}$ for GD on nonconvex functions.

Similarly, if we perform SGD update $x_{t+1}\leftarrow x_t - \eta\nabla f_i(x_t)$ each time for a random $i\in[n]$, then

$$\begin{aligned}
f(x_t) - \mathbb{E}_i[f(x_{t+1})] &\geq \mathbb{E}_i\big[\langle \nabla f(x_t), x_t - x_{t+1}\rangle - \frac{L}{2}\|x_t - x_{t+1}\|^2\big] = \eta\|\nabla f(x_t)\|^2 - \frac{\eta^2 L}{2}\mathbb{E}_i\big[\|\nabla f_i(x_t)\|^2\big] \\
&= \big(\eta - \frac{\eta^2 L}{2}\big)\|\nabla f(x_t)\|^2 - \frac{\eta^2 L}{2}\mathbb{E}_i\big[\|\nabla f_i(x_t) - \nabla f(x_t)\|^2\big] \\
&\geq \big(\eta - \frac{\eta^2 L}{2}\big)\|\nabla f(x_t)\|^2 - \frac{\eta^2 L}{2}\mathcal{V} \ .
\end{aligned}$$

Therefore, choosing $\eta = \min\big\{\frac{1}{L}, \frac{\varepsilon^2}{L\mathcal{V}}\big\}$, we can conclude that if $t$ is randomly chosen among $t = 0,1,2,\dots,T-1$, then it satisfies $\mathbb{E}[\|\nabla f(x_t)\|^2] \leq \varepsilon^2$ if $T \geq \Omega\big(\big(\frac{L}{\varepsilon^2} + \frac{L\mathcal{V}}{\varepsilon^4}\big)(f(x_0) - f(x^*))\big)$. This is the $T \propto \frac{1}{\varepsilon^4}$ convergence rate for SGD. One can use acceleration techniques to improve the lower-order term $\varepsilon^{-2}$ in this complexity [23], but not the $\varepsilon^{-4}$ term.

# C Missing Proofs for Section 5: Auxiliary Lemmas

## C.1 Proof of Lemma 5.3

**Lemma 5.3.** *There exists absolute constant $C > 1$ such that for any $x \in \mathbb{R}^d$, $T \geq 1$, $p \in (0,1)$:*

- *if we run Oja's algorithm once for $T$ iterations, with $\eta = \Theta(\sqrt{T})$, we can find unit vector $y$ such that, with at with probability at least $4/5$,*

$$y^\top \nabla^2 f(x) y \le \lambda_{\min}(\nabla^2 f(x)) + C \cdot \tfrac{L \log(d)}{\sqrt{T}} \enspace .$$

- *if we run Oja's algorithm $O(\log(1/p))$ times, then with probability at least $1-p$, we can*

$$\text{either conclude} \qquad \lambda_{\min}(\nabla^2 f(x)) \ge -C \cdot \frac{L \log(d/p)}{\sqrt{T}} \enspace ,$$

$$\text{or find } y \in \mathbb{R}^d \text{ such that} \qquad y^\top \nabla^2 f(x) y \le -\frac{C}{2} \cdot \frac{L \log(d/p)}{\sqrt{T}} \enspace .$$

*The total number of hessian-vector products is at most $O(T \log(1/p))$.*

*Proof of Lemma 5.3.* It is clear that all matrices generated from $\mathcal{D}$ are symmetric, and have eigenvalues between 0 and 1. By applying Theorem 5.1, and setting $y$ to be a uniform random one among $w_1, \dots, w_T$, we have with probability at least $9/10$ (over the randomness of $\mathcal{D}$):

$$\mathbb{E}_y\big[\lambda_{\max}(\mathbf{B}) - y^\top \mathbf{B} y\big] \le O\big(\log(d)/\sqrt{T}\big)$$

By Markov's bound (and noting that $\lambda_{\max}(\mathbf{B}) - y^\top \mathbf{B} y$ is always non-negative), with probability at least $4/5$, we have $\lambda_{\max}(\mathbf{B}) - y^\top \mathbf{B} y \le O\big(\log(d)/\sqrt{T}\big)$. This finishes the proof of the first item, after plugging in the definition of $\mathbf{B} = \frac{L \cdot I - \nabla^2 f(x)}{2L}$.

For the second item, suppose we run Oja's algorithm, independently, for $O(\log(1/p))$ times, and let the output vector $y$ be denoted as $y_t$ for each run $t \in [O(\log(1/p))]$. We know that, with probability at least $1 - p/2$, at least one of the runs is successful and outputs $y_t$ satisfying $\lambda_{\max}(\mathbf{B}) - y_t^\top \mathbf{B} y_t \le O\big(\log(d)/\sqrt{T}\big)$.

Moreover, to test whether the $t$-th run is successful, we generate additionally $T$ copies of samples from $\mathcal{D}$, denoted by $\mathbf{B}_{t,1}, \dots, \mathbf{B}_{t,T}$ from $\mathcal{D}$. By Bernstein's inequality, we have for every $\varepsilon \in (0,1)$:

$$\Pr_{\mathbf{B}_{t,1},\dots,\mathbf{B}_{t,T}}\left[\left| y_t^\top \frac{\mathbf{B}_{t,1} + \cdots + \mathbf{B}_{t,T}}{T} y_t - y_t^\top \mathbf{B} y_t \right| > \varepsilon\right] \le e^{-\Omega(T\varepsilon^2)} \enspace .$$

In other words, by union bound, with probability at least $1 - p/2$, we have

$$\left| y_t^\top \frac{\mathbf{B}_{t,1} + \cdots + \mathbf{B}_{t,T}}{T} y_t - y_t^\top \mathbf{B} y_t \right| \le O\big(\frac{\log(1/p)}{\sqrt{T}}\big) \qquad \forall t \in [O(\log(1/p))] \enspace .$$

Conditioning on that both two events hold (with probability $\ge 1-p$), define

$$t^* = \arg\min_t \{ y_t^\top \frac{\mathbf{B}_{t,1} + \cdots + \mathbf{B}_{t,T}}{T} y_t \}, \qquad y \stackrel{\text{def}}{=} y_{t^*}, \qquad \text{and } \beta \stackrel{\text{def}}{=} y_{t^*}^\top \frac{\mathbf{B}_{t^*,1} + \cdots + \mathbf{B}_{t^*,T}}{T} y_{t^*}.$$

We conclude that, there exist some constant $C > 1$ such that

$$y^\top \mathbf{B} y \ge \lambda_{\max}(\mathbf{B}) - O\big(\log d/\sqrt{T}\big) - O\big(\log(1/p)/\sqrt{T}\big) \ge \lambda_{\max}(\mathbf{B}) - C \cdot \frac{\log(d/p)}{\sqrt{T}}$$

$$|\beta - y^\top \mathbf{B} y| \le O\big(\log(1/p)/\sqrt{T}\big) \le C \cdot \frac{\log(d/p)}{\sqrt{T}}$$

Plugging in the definition $\mathbf{B} = \frac{L \cdot I - \nabla^2 f(x)}{2L}$ and choosing $\rho = L - 2L\beta$, we have

$$y^\top \nabla^2 f(x) y \le \lambda_{\min}(\nabla^2 f(x)) + C \cdot \frac{2L \log(d/p)}{\sqrt{T}}$$

$$|\rho - y^\top \nabla^2 f(x) y| \le C \cdot \frac{2L \log(d/p)}{\sqrt{T}}$$

Finally, if $\rho < -4C \cdot \frac{2L\log(d/p)}{\sqrt{T}}$, then we have $y^\top \nabla^2 f(x) y \leq -3C \cdot \frac{2L\log(d/p)}{\sqrt{T}}$; otherwise, if $\rho \geq -4C \cdot \frac{2L\log(d/p)}{\sqrt{T}}$, then we conclude that $\lambda_{\min}(\nabla^2 f(x)) \geq -6C \cdot \frac{2L\log(d/p)}{\sqrt{T}}$. This finishes the proof of the second item. □

## C.2 Proof of Lemma 5.5

The following definition and properties of Fenchel dual for convex functions is classical, and can be found for instance in the textbook [43].

**Definition C.1.** *Given proper convex function $h(y)$, its Fenchel dual $h^*(\beta) \stackrel{\text{def}}{=} \max_y\{y^\top\beta - h(y)\}$.*

**Proposition C.2.** $\nabla h^*(\beta) = \arg\max_y\{y^\top\beta - h(y)\}$.

**Proposition C.3.** *If $h(\cdot)$ is $\sigma$-strongly convex, then $h^*(\cdot)$ is $\frac{1}{\sigma}$-smooth.*

**Lemma 5.5.** *Let $x^*$ be the unique minimizer of $G(y)$, and $x$ be an arbitrary vector in the domain of $\{x \in \mathbb{R}^d : \psi(x) < +\infty\}$. Then, for every $\eta \in \big(0, \frac{1}{L+2\sigma}\big]$, we have*

$$\|\mathcal{G}_{F,\eta}(x)\|^2 + \sigma^2\|x - \widehat{\mathsf{x}}\|^2 \leq O\big(\sigma^2\|x^* - \widehat{\mathsf{x}}\|^2 + \|\mathcal{G}_{G,\eta}(x)\|^2\big) \ .$$

*Proof of Lemma 5.5.* Define

$$z = \arg\min_y \big\{\psi(y) + \langle\nabla f(x), y\rangle + \frac{1}{2\eta}\|y - x\|^2\big\}$$
$$\overline{z} = \arg\min_y \big\{\psi(y) + \langle\nabla f(x) + 2\sigma(x - \widehat{\mathsf{x}}), y\rangle + \frac{1}{2\eta}\|y - x\|^2\big\}$$

We have by definition $\mathcal{G}_{F,\eta}(x) = \frac{x-z}{\eta}$ and $\mathcal{G}_{G,\eta}(x) = \frac{x-\overline{z}}{\eta}$. Therefore, by AM-GM,

$$\|\mathcal{G}_{F,\eta}(x)\|^2 \leq 2\|\mathcal{G}_{G,\eta}(x)\|^2 + \frac{2}{\eta^2}\|z - \overline{z}\|^2 \ . \tag{C.1}$$

On the other hand, let us denote by $h(y) \stackrel{\text{def}}{=} \psi(y) + \frac{1}{2\eta}\|y\|^2$ and recall the definition of Fenchel dual $h^*(\beta) = \max_y\{y^\top\beta - h(y)\}$. Proposition C.2 says $\nabla h^*(\beta) = \max_y\{y^\top\beta - h(y)\}$. This implies

$$z = \nabla h^*\big(\tfrac{x}{\eta} - \nabla f(x)\big) \quad \text{and} \quad \overline{z} = \nabla h^*\big(\tfrac{x}{\eta} - \nabla f(x) - 2\sigma(x - \widehat{\mathsf{x}})\big) \ .$$

Using the property that $h^*(\cdot)$ is $\eta$-smooth (because $h(y)$ is $1/\eta$-strongly convex, see Proposition C.3), we have

$$\frac{1}{\eta^2}\|z - \overline{z}\|^2 \leq \|2\sigma(x - \widehat{\mathsf{x}})\|^2 \leq 8\sigma^2\|x^* - \widehat{\mathsf{x}}\|^2 + 8\sigma^2\|x - x^*\|^2 \ . \tag{C.2}$$

Next, recall the following property about gradient mapping —see for instance [47, Lemma 3.7])— for every $x^*$:[22]

$$\forall \eta \leq \frac{1}{L + 2\sigma}: \quad G(x^*) \geq G(\overline{z}) + \langle\mathcal{G}_{G,\eta}(x), x^* - x\rangle + \frac{\eta}{2}\|\mathcal{G}_{G,\eta}(x)\|^2 + \frac{\sigma}{2}\|x^* - x\|^2 \ .$$

Using $G(x^*) \leq G(\overline{z})$, the non-negativity of $\|\mathcal{G}_{G,\eta}(x)\|^2$, and Young's inequality $|\langle\mathcal{G}_{G,\eta}(x), x^* - x\rangle| \leq \frac{1}{\sigma}\|\mathcal{G}_{G,\eta}(x)\|^2 + \frac{\sigma}{4}\|x - x^*\|^2$, we have

$$\frac{\sigma^2}{4}\|x - x^*\|^2 \leq \|\mathcal{G}_{G,\eta}(x)\|^2 \ . \tag{C.3}$$

Finally, combining (C.1), (C.2), and (C.3), we have the desired result. □

[22] To apply [47, Lemma 3.7], we observe that $g(x) = f(s) + \sigma\|x - \widehat{\mathsf{x}}\|^2$ is convex and $(L + 2\sigma)$-smooth.

# D Missing Proofs for Section 6: Natasha 1.5

## D.1 Proof of Lemma 6.1

**Lemma 6.1.** *We have* $\mathbb{E}_{i_t^s}\big[\|\widetilde{\nabla} f^s(x_t^s) - \nabla f^s(x_t^s)\|^2\big] \leq pL^2\|x_t^s - \widehat{\mathsf{x}}^s\|^2 + pL^2\sum_{k=0}^{s-1}\|\widehat{\mathsf{x}}^k - \widehat{\mathsf{x}}^{k+1}\|^2 + \|\mathbf{e}\|^2$ .

*Proof.* We have

$$\mathbb{E}_{i_t^s}\big[\|\widetilde{\nabla} f^s(x_t^s) - \nabla f^s(x_t^s)\|^2\big] = \mathbb{E}_{i_t^s}\big[\|\widetilde{\nabla} f(x_t^s) - \nabla f(x_t^s)\|^2\big]$$

$$= \mathbb{E}_{i\in_R[n]}\Big[\big\|\big(\nabla f_i(x_t^s) - \nabla f_i(\widetilde{\mathbf{x}})\big) - \big(\nabla f(x_t^s) - \nabla f(\widetilde{\mathbf{x}})\big) + \mathbf{e}\big\|^2\Big]$$

$$\overset{①}{=} \mathbb{E}_{i\in_R[n]}\Big[\big\|\big(\nabla f_i(x_t^s) - \nabla f_i(\widetilde{\mathbf{x}})\big) - \big(\nabla f(x_t^s) - \nabla f(\widetilde{\mathbf{x}})\big)\big\|^2\Big] + \|\mathbf{e}\|^2$$

$$\overset{②}{\leq} \mathbb{E}_{i\in_R[n]}\Big[\big\|\nabla f_i(x_t^s) - \nabla f_i(\widetilde{\mathbf{x}})\big\|^2\Big] + \|\mathbf{e}\|^2$$

$$\overset{③}{\leq} p\mathbb{E}_{i\in_R[n]}\Big[\big\|\nabla f_i(x_t^s) - \nabla f_i(\widehat{\mathsf{x}}^s)\big\|^2\Big] + p\sum_{k=0}^{s-1}\mathbb{E}_{i\in_R[n]}\Big[\big\|\nabla f_i(\widehat{\mathsf{x}}^k) - \nabla f_i(\widehat{\mathsf{x}}^{k+1})\big\|^2\Big] + \|\mathbf{e}\|^2$$

$$\overset{④}{\leq} pL^2\|x_t^s - \widehat{\mathsf{x}}^s\|^2 + pL^2\sum_{k=0}^{s-1}\|\widehat{\mathsf{x}}^k - \widehat{\mathsf{x}}^{k+1}\|^2 + \|\mathbf{e}\|^2 \ .$$

Above, equality ① is because $\mathbb{E}[\|a+b\|^2] = \mathbb{E}[\|a\|^2] + \|b\|^2$ for any random vector $a$ and non-random vector $b$, as long as $\mathbb{E}[a] = \vec{0}$; inequality ② is because for any random vector $\zeta \in \mathbb{R}^d$, it holds that $\mathbb{E}\|\zeta - \mathbb{E}\zeta\|^2 = \mathbb{E}\|\zeta\|^2 - \|\mathbb{E}\zeta\|^2$; inequality ③ is because $\widehat{\mathsf{x}}^0 = \widetilde{\mathbf{x}}$ and for any $p$ vectors $a_1, a_2, \dots, a_p \in \mathbb{R}^d$, it holds that $\|a_1 + \dots + a_p\|^2 \leq p\|a_1\|^2 + \dots + p\|a_p\|^2$; and inequality ④ is because each $f_i(\cdot)$ is $L$-smooth. □

## D.2 Proof of Lemma 6.4

The following inequality is classically known as the "regret inequality" for proximal mirror descent, and its proof is classical.

**Fact D.1.** *If* $x_{t+1} = \arg\min_{y\in\mathbb{R}^d}\{\psi(y) + \frac{1}{2\alpha}\|y - x_t\|^2 + \langle w, y\rangle\}$*, then for every* $u \in \mathbb{R}^d$*:*

$$\langle w, x_{t+1} - u\rangle + \psi(x_{t+1}) - \psi(u) \leq \frac{\|x_t - u\|^2}{2\alpha} - \frac{\|x_{t+1} - u\|^2}{2\alpha} - \frac{\|x_{t+1} - x_t\|^2}{2\alpha} \ .$$

*Proof.* Recall that the minimality of $x_{t+1} = \arg\min_{y\in\mathbb{R}^d}\{\frac{1}{2\alpha}\|y - x_t\|^2 + \psi(y) + \langle w, y\rangle\}$ implies the existence of some subgradient $g \in \partial\psi(x_{t+1})$ which satisfies $\frac{1}{\alpha}(x_{t+1} - x_t) + w + g = 0$. Combining this with $\psi(u) - \psi(x_{t+1}) \geq \langle g, u - x_{t+1}\rangle$, which is due to the convexity of $\psi(\cdot)$, we immediately have $\psi(u) - \psi(x_{t+1}) + \langle \frac{1}{\alpha}(x_{t+1} - x_t) + w, u - x_{t+1}\rangle \geq \langle \frac{1}{\alpha}(x_{t+1} - x_t) + w + g, u - x_{t+1}\rangle = 0$. Rearranging this inequality we have

$$\begin{aligned}\langle w, x_{t+1} - u\rangle + \psi(x_{t+1}) - \psi(u) &\leq \langle -\frac{1}{\alpha}(x_{t+1} - x_t), x_{t+1} - u\rangle \\ &= \frac{\|x_t - u\|^2}{2\alpha} - \frac{\|x_{t+1} - u\|^2}{2\alpha} - \frac{\|x_{t+1} - x_t\|^2}{2\alpha} \ .\end{aligned}$$

□

**Lemma 6.4.** *As long as* $\alpha \leq \frac{1}{2L+4\sigma}$*, letting* $x_*^s = \arg\min_x\{F(x) + \sigma\|x - \widehat{\mathsf{x}}^s\|^2\}$*, we have*

$$\mathbb{E}\Big[\big(F^s(\widehat{\mathsf{x}}^{s+1}) - F^s(x_*^s)\big)\Big] \leq \mathbb{E}\Big[\frac{F^s(\widehat{\mathsf{x}}^s) - F^s(x_*^s)}{\sigma\alpha m/4} + 2\alpha pL^2\Big(\sum_{k=0}^{s}\|\widehat{\mathsf{x}}^k - \widehat{\mathsf{x}}^{k+1}\|^2\Big)\Big] + \frac{3}{\sigma}\|\mathbf{e}\|^2 \ .$$

*Proof of Lemma 6.4.* We first compute that

$$
\begin{aligned}
&F^s(x_{t+1}^s) - F^s(u) = f^s(x_{t+1}^s) - f^s(u) + \psi(x_{t+1}^s) - \psi(u) \\
&\overset{①}{\leq} f^s(x_t^s) + \langle \nabla f^s(x_t^s), x_{t+1}^s - x_t^s \rangle + \frac{L+2\sigma}{2}\|x_t^s - x_{t+1}^s\|^2 - f^s(u) + \psi(x_{t+1}^s) - \psi(u) \\
&\overset{②}{\leq} \langle \nabla f^s(x_t^s), x_{t+1}^s - x_t^s \rangle + \frac{L+2\sigma}{2}\|x_t^s - x_{t+1}^s\|^2 + \langle \nabla f^s(x_t^s), x_t^s - u\rangle + \psi(x_{t+1}^s) - \psi(u) \ . \qquad \text{(D.1)}
\end{aligned}
$$

Above, inequality ① uses the fact that $f^s(\cdot)$ is $(L+2\sigma)$-smooth; and inequality ② uses the convexity of $f^s(\cdot)$. Now, we take expectation with respect to $i_t^s$ on both sides of (D.1), and derive that:

$$
\begin{aligned}
&\mathbb{E}_{i_t^s}\big[F^s(x_{t+1}^s)\big] - F^s(u) + \langle \mathbf{e}, x_t^s - u\rangle \\
&\overset{①}{\leq} \mathbb{E}_{i_t^s}\Big[\langle \widetilde{\nabla} f^s(x_t^s) - \nabla f^s(x_t^s), x_t^s - x_{t+1}^s\rangle + \langle \widetilde{\nabla} f^s(x_t^s), x_{t+1}^s - u\rangle + \frac{L+2\sigma}{2}\|x_t^s - x_{t+1}^s\|^2 + \psi(x_{t+1}^s) - \psi(u)\Big] \\
&\overset{②}{\leq} \mathbb{E}_{i_t^s}\Big[\langle \widetilde{\nabla} f^s(x_t^s) - \nabla f^s(x_t^s), x_t^s - x_{t+1}^s\rangle + \frac{\|x_t^s - u\|^2}{2\alpha} - \frac{\|x_{t+1}^s - u\|^2}{2\alpha} - \big(\frac{1}{2\alpha} - \frac{L+2\sigma}{2}\big)\|x_{t+1}^s - x_t^s\|^2\Big] \\
&\overset{③}{\leq} \mathbb{E}_{i_t^s}\Big[\alpha\big\|\widetilde{\nabla} f^s(x_t^s) - \nabla f^s(x_t^s)\big\|^2 + \frac{\|x_t^s - u\|^2}{2\alpha} - \frac{\|x_{t+1}^s - u\|^2}{2\alpha}\Big] \\
&\overset{④}{\leq} \mathbb{E}_{i_t^s}\Big[\alpha p L^2\|x_t^s - \widehat{\mathsf{x}}^s\|^2 + \alpha p L^2 \sum_{k=0}^{s-1}\|\widehat{\mathsf{x}}^k - \widehat{\mathsf{x}}^{k+1}\|^2 + \alpha\|\mathbf{e}\|^2 + \frac{\|x_t^s - u\|^2}{2\alpha} - \frac{\|x_{t+1}^s - u\|^2}{2\alpha}\Big] \ . \qquad \text{(D.2)}
\end{aligned}
$$

Above, inequality ① follows from (D.1) together with the following inequality (noticing that $x_t^s$ and $u$ do not depend on the randomness of $i_t^s$, and $\mathbb{E}_{i_t^s}[\widetilde{\nabla} f^s(x_t^s)] = \nabla f^s(x_t^s) + \mathbf{e}$):

$$
\begin{aligned}
&\mathbb{E}_{i_t^s}\big[\langle \nabla f^s(x_t^s), x_{t+1}^s - x_t^s\rangle + \langle \nabla f^s(x_t^s), x_t^s - u\rangle\big] \\
&\qquad = \mathbb{E}_{i_t^s}\big[\langle \widetilde{\nabla} f^s(x_t^s) - \nabla f^s(x_t^s), x_t^s - x_{t+1}^s\rangle + \langle \widetilde{\nabla} f^s(x_t^s), x_{t+1}^s - u\rangle\big] - \langle \mathbf{e}, x_t^s - u\rangle \ ;
\end{aligned}
$$

inequality ② uses Fact D.1; inequality ③ uses $\alpha \leq \frac{1}{2L+4\sigma}$ together with Young's inequality $\langle a,b\rangle \leq \frac{1}{2}\|a\|^2 + \frac{1}{2}\|b\|^2$; and inequality ④ uses Lemma 6.1.

Next, choosing $u = x_*^s$ to be the unique minimizer of $F^s(\cdot) = f^s(\cdot) + \psi(\cdot)$, and telescoping inequality (D.2) for $t = 0, 1, \dots, m-1$, we have (noticing that $F^s(x_m^s) - F^s(x_*^s) \geq 0$)

$$
\begin{aligned}
&\mathbb{E}\Big[\sum_{t=1}^{m-1}\big(F^s(x_t^s) - F^s(x_*^s)\big) + \sum_{t=0}^{m-1}\langle \mathbf{e}, x_t^s - x_*^s\rangle\Big] - \alpha m\|\mathbf{e}\|^2 \\
&\leq \mathbb{E}\Big[\frac{\|x_0^s - x_*^s\|^2}{2\alpha} + \sum_{t=0}^{m-1}\Big(\alpha p L^2\|x_t^s - \widehat{\mathsf{x}}^s\|^2 + \alpha p L^2\sum_{k=0}^{s-1}\|\widehat{\mathsf{x}}^k - \widehat{\mathsf{x}}^{k+1}\|^2\Big)\Big] \\
&\leq \mathbb{E}\Big[\frac{F^s(\widehat{\mathsf{x}}^s) - F^s(x_*^s)}{\sigma\alpha} + \alpha p m L^2\Big(\sum_{k=0}^{s}\|\widehat{\mathsf{x}}^k - \widehat{\mathsf{x}}^{k+1}\|^2\Big)\Big] \ . \qquad \text{(D.3)}
\end{aligned}
$$

Above, the second inequality uses the fact that $\widehat{\mathsf{x}}^{s+1}$ is chosen from $\{x_0^s, \dots, x_{m-1}^s\}$ uniformly at random, as well as the $\sigma$-strong convexity of $F^s(\cdot)$.

At this point, we apply Young's inequality

$$
-\langle \mathbf{e}, x_t^s - x_*^s\rangle \leq \frac{\sigma\|x_t^s - x_*^s\|^2}{4} + \frac{1}{\sigma}\|\mathbf{e}\|^2 \leq \frac{1}{2}\big(F^s(x_t^s) - F^s(x_*^s)\big) + \frac{1}{\sigma}\|\mathbf{e}\|^2 \qquad \text{(D.4)}
$$

and are ready to plug (D.4) into (D.3). Using again the fact that that $\widehat{\mathsf{x}}^{s+1}$ is chosen from $\{x_0^s, \dots, x_{m-1}^s\}$ uniformly at random, we have

$$\mathbb{E}\Big[\frac{1}{2}m\big(F^s(\widehat{\mathsf{x}}^{s+1}) - F^s(x^s_*)\big)\Big] - \big(\alpha m + \frac{m}{\sigma}\big)\|\mathbf{e}\|^2$$
$$\leq F^s(x^s_0) - F^s(x^s_*) + \mathbb{E}\Big[\frac{F^s(\widehat{\mathsf{x}}^s) - F^s(x^s_*)}{\sigma\alpha} + \alpha p m L^2\Big(\sum_{k=0}^{s}\|\widehat{\mathsf{x}}^k - \widehat{\mathsf{x}}^{k+1}\|^2\Big)\Big] \ . \tag{D.5}$$

Finally, using our choices $x^s_0 = \widehat{\mathsf{x}}^s$ and $\frac{1}{2\sigma\alpha} \geq 1$ (which is implied by $\alpha \leq \frac{1}{2L+4\sigma}$), we divide both sides of (D.5) by $m/2$, and rearrange the terms:

$$\mathbb{E}\Big[\big(F^s(\widehat{\mathsf{x}}^{s+1}) - F^s(x^s_*)\big)\Big] \leq \mathbb{E}\Big[\frac{F^s(\widehat{\mathsf{x}}^s) - F^s(x^s_*)}{\sigma\alpha m/4} + 2\alpha p L^2\Big(\sum_{k=0}^{s}\|\widehat{\mathsf{x}}^k - \widehat{\mathsf{x}}^{k+1}\|^2\Big)\Big] + \frac{3}{\sigma}\|\mathbf{e}\|^2 \ . \qquad \square$$

## D.3 Proof of Lemma 6.5

**Lemma 6.5.** *If $\alpha \leq \frac{1}{2L+4\sigma}$, $\alpha \geq \frac{8}{\sigma m}$ and $\alpha \leq \frac{\sigma}{4p^2L^2}$, we have*

$$\sum_{s=0}^{p-1}\mathbb{E}\Big[\sigma\|\widehat{\mathsf{x}}^s - \widehat{\mathsf{x}}^{s+1}\|^2 + \frac{\sigma}{2}\|\widehat{\mathsf{x}}^s - x^s_*\|^2\Big] \leq 2\mathbb{E}\Big[F(\widehat{\mathsf{x}}^0) - F(\widehat{\mathsf{x}}^p)\Big] + \frac{3p\mathcal{V}}{\sigma B} \ ,$$

*where recall $x^s_* \stackrel{\text{def}}{=} \arg\min_x\{F(x) + \sigma\|x - \widehat{\mathsf{x}}^s\|^2\}$.*

*Proof of Lemma 6.5.* Telescoping Lemma 6.4 for all the subepochs $s = 0, 1, \dots, p-1$, we have

$$\sum_{s=0}^{p-1}\mathbb{E}\Big[\frac{\sigma}{2}\|\widehat{\mathsf{x}}^s - \widehat{\mathsf{x}}^{s+1}\|^2 + \big(F^s(\widehat{\mathsf{x}}^{s+1}) - F^s(x^s_*)\big)\Big]$$
$$\leq \sum_{s=0}^{p-1}\mathbb{E}\Big[\frac{\sigma}{2}\|\widehat{\mathsf{x}}^s - \widehat{\mathsf{x}}^{s+1}\|^2 + \frac{F^s(\widehat{\mathsf{x}}^s) - F^s(x^s_*)}{\sigma\alpha m/4} + 2\alpha p^2L^2\|\widehat{\mathsf{x}}^s - \widehat{\mathsf{x}}^{s+1}\|^2 + \frac{3p}{\sigma}\|\mathbf{e}\|^2\Big]$$
$$\overset{①}{\leq} \sum_{s=0}^{p-1}\mathbb{E}\Big[\frac{F^s(\widehat{\mathsf{x}}^s) - F^s(x^s_*)}{\sigma\alpha m/4} + \sigma\cdot\|\widehat{\mathsf{x}}^{s+1} - \widehat{\mathsf{x}}^s\|^2 + \frac{3p}{\sigma}\|\mathbf{e}\|^2\Big]$$
$$\overset{②}{=} \sum_{s=0}^{p-1}\mathbb{E}\Big[\frac{F^s(\widehat{\mathsf{x}}^s) - F^s(x^s_*)}{\sigma\alpha m/4} + \big(F^s(\widehat{\mathsf{x}}^{s+1}) - F^s(\widehat{\mathsf{x}}^s)\big) - \big(F(\widehat{\mathsf{x}}^{s+1}) - F(\widehat{\mathsf{x}}^s)\big) + \frac{3p}{\sigma}\|\mathbf{e}\|^2\Big]$$

Above, ① uses $4\alpha p^2L^2 \leq \sigma$, and ② uses the definition $F^s(y) = F(y) + \sigma\|y - \widehat{\mathsf{x}}^s\|^2$.

Finally, rearranging both sides, and using the fact that $\frac{1}{\sigma\alpha m} \leq \frac{1}{8}$ and the fact that $\mathbb{E}[\|\mathbf{e}\|^2] \leq \frac{\mathcal{V}}{B}$ from Claim 6.2, we have

$$\sum_{s=0}^{p-1}\mathbb{E}\Big[\sigma\|\widehat{\mathsf{x}}^s - \widehat{\mathsf{x}}^{s+1}\|^2 + \big(F^s(\widehat{\mathsf{x}}^s) - F^s(x^s_*)\big)\Big] \leq 2\mathbb{E}\Big[F(\widehat{\mathsf{x}}^0) - F(\widehat{\mathsf{x}}^p)\Big] + \frac{3p\mathcal{V}}{\sigma B} \ .$$

If we further apply the $\sigma$-strong convexity of $F^s(\cdot)$ we have the desired inequality. $\square$

# References

[1] Naman Agarwal, Zeyuan Allen-Zhu, Brian Bullins, Elad Hazan, and Tengyu Ma. Finding Approximate Local Minima for Nonconvex Optimization in Linear Time. In *STOC*, 2017. Full version available at `http://arxiv.org/abs/1611.01146`.

[2] Zeyuan Allen-Zhu. Katyusha: The First Direct Acceleration of Stochastic Gradient Methods. In *STOC*, 2017. Full version available at `http://arxiv.org/abs/1603.05953`.

[3] Zeyuan Allen-Zhu. Natasha: Faster Non-Convex Stochastic Optimization via Strongly Non-Convex Parameter. In *ICML*, 2017. Full version available at http://arxiv.org/abs/1702.00763.

[4] Zeyuan Allen-Zhu. Katyusha X: Practical Momentum Method for Stochastic Sum-of-Nonconvex Optimization. In *ICML*, 2018. Full version available at http://arxiv.org/abs/1802.03866.

[5] Zeyuan Allen-Zhu. How To Make the Gradients Small Stochastically: Even Faster Convex and Nonconvex SGD. *ArXiv e-prints*, abs/1801.02982, January 2018. Full version available at http://arxiv.org/abs/1801.02982.

[6] Zeyuan Allen-Zhu and Elad Hazan. Variance Reduction for Faster Non-Convex Optimization. In *ICML*, 2016. Full version available at http://arxiv.org/abs/1603.05643.

[7] Zeyuan Allen-Zhu and Yuanzhi Li. LazySVD: Even Faster SVD Decomposition Yet Without Agonizing Pain. In *NIPS*, 2016. Full version available at http://arxiv.org/abs/1607.03463.

[8] Zeyuan Allen-Zhu and Yuanzhi Li. First Efficient Convergence for Streaming k-PCA: a Global, Gap-Free, and Near-Optimal Rate. In *FOCS*, 2017. Full version available at http://arxiv.org/abs/1607.07837.

[9] Zeyuan Allen-Zhu and Yuanzhi Li. Follow the Compressed Leader: Faster Online Learning of Eigenvectors and Faster MMWU. In *ICML*, 2017. Full version available at http://arxiv.org/abs/1701.01722.

[10] Zeyuan Allen-Zhu and Yuanzhi Li. Neon2: Finding Local Minima via First-Order Oracles. *ArXiv e-prints*, abs/1711.06673, November 2017. Full version available at http://arxiv.org/abs/1711.06673.

[11] Zeyuan Allen-Zhu and Lorenzo Orecchia. Linear Coupling: An Ultimate Unification of Gradient and Mirror Descent. In *Proceedings of the 8th Innovations in Theoretical Computer Science*, ITCS '17, 2017. Full version available at http://arxiv.org/abs/1407.1537.

[12] Sanjeev Arora, Rong Ge, Tengyu Ma, and Ankur Moitra. Simple, Efficient, and Neural Algorithms for Sparse Coding. In *COLT*, 2015.

[13] Yair Carmon, John C. Duchi, Oliver Hinder, and Aaron Sidford. Accelerated Methods for Non-Convex Optimization. *ArXiv e-prints*, abs/1611.00756, November 2016.

[14] Yair Carmon, Oliver Hinder, John C. Duchi, and Aaron Sidford. "Convex Until Proven Guilty": Dimension-Free Acceleration of Gradient Descent on Non-Convex Functions. In *ICML*, 2017.

[15] Yuxin Chen and Emmanuel Candes. Solving random quadratic systems of equations is nearly as easy as solving linear systems. In *Advances in Neural Information Processing Systems*, pages 739–747, 2015.

[16] Anna Choromanska, Mikael Henaff, Michael Mathieu, Gérard Ben Arous, and Yann LeCun. The loss surfaces of multilayer networks. In *AISTATS*, 2015.

[17] Yann N Dauphin, Razvan Pascanu, Caglar Gulcehre, Kyunghyun Cho, Surya Ganguli, and Yoshua Bengio. Identifying and attacking the saddle point problem in high-dimensional non-convex optimization. In *NIPS*, pages 2933–2941, 2014.

[18] Aaron Defazio, Francis Bach, and Simon Lacoste-Julien. SAGA: A Fast Incremental Gradient Method With Support for Non-Strongly Convex Composite Objectives. In *NIPS*, 2014.

[19] John Duchi, Elad Hazan, and Yoram Singer. Adaptive subgradient methods for online learning and stochastic optimization. *The Journal of Machine Learning Research*, 12:2121–2159, 2011.

[20] Roy Frostig, Rong Ge, Sham M. Kakade, and Aaron Sidford. Un-regularizing: approximate proximal point and faster stochastic algorithms for empirical risk minimization. In *ICML*, 2015.

[21] Dan Garber, Elad Hazan, Chi Jin, Sham M. Kakade, Cameron Musco, Praneeth Netrapalli, and Aaron Sidford. Robust shift-and-invert preconditioning: Faster and more sample efficient algorithms for eigenvector computation. In *ICML*, 2016.

[22] Rong Ge, Furong Huang, Chi Jin, and Yang Yuan. Escaping from saddle points—online stochastic gradient for tensor decomposition. In *Proceedings of the 28th Annual Conference on Learning Theory*, COLT 2015, 2015.

[23] Saeed Ghadimi and Guanghui Lan. Accelerated gradient methods for nonconvex nonlinear and stochastic programming. *Mathematical Programming*, pages 1–26, feb 2015. ISSN 0025-5610.

[24] I. J. Goodfellow, O. Vinyals, and A. M. Saxe. Qualitatively characterizing neural network optimization

problems. *ArXiv e-prints*, December 2014.

[25] Elad Hazan, Kfir Yehuda Levy, and Shai Shalev-Shwartz. On graduated optimization for stochastic non-convex problems. In *International Conference on Machine Learning*, pages 1833–1841, 2016.

[26] Xi He, Dheevatsa Mudigere, Mikhail Smelyanskiy, and Martin Takáč. Distributed Hessian-Free Optimization for Deep Neural Network. *ArXiv e-prints*, abs/1606.00511, June 2016.

[27] Chi Jin, Rong Ge, Praneeth Netrapalli, Sham M Kakade, and Michael I Jordan. How to Escape Saddle Points Efficiently. In *ICML*, 2017.

[28] Rie Johnson and Tong Zhang. Accelerating stochastic gradient descent using predictive variance reduction. In *Advances in Neural Information Processing Systems*, NIPS 2013, pages 315–323, 2013.

[29] Diederik Kingma and Jimmy Ba. Adam: A method for stochastic optimization. *ArXiv e-prints*, abs/1412.6980, 12 2014.

[30] Lihua Lei, Cheng Ju, Jianbo Chen, and Michael I Jordan. Nonconvex Finite-Sum Optimization Via SCSG Methods. In *NIPS*, 2017.

[31] Yuanzhi Li and Yang Yuan. Convergence Analysis of Two-layer Neural Networks with ReLU Activation. In *NIPS*, 2017.

[32] Hongzhou Lin, Julien Mairal, and Zaid Harchaoui. A Universal Catalyst for First-Order Optimization. In *NIPS*, 2015.

[33] Yurii Nesterov. *Introductory Lectures on Convex Programming Volume: A Basic course*, volume I. Kluwer Academic Publishers, 2004. ISBN 1402075537.

[34] Yurii Nesterov. Accelerating the cubic regularization of newton's method on convex problems. *Mathematical Programming*, 112(1):159–181, 2008.

[35] Yurii Nesterov. How to make the gradients small. *Optima*, 88:10–11, 2012.

[36] Yurii Nesterov and Boris T. Polyak. Cubic regularization of newton method and its global performance. *Mathematical Programming*, 108(1):177–205, 2006.

[37] Erkki Oja. Simplified neuron model as a principal component analyzer. *Journal of mathematical biology*, 15(3):267–273, 1982.

[38] Barak A Pearlmutter. Fast exact multiplication by the hessian. *Neural computation*, 6(1):147–160, 1994.

[39] Sashank J. Reddi, Ahmed Hefny, Suvrit Sra, Barnabas Poczos, and Alex Smola. Stochastic variance reduction for nonconvex optimization. In *ICML*, 2016.

[40] Sashank J Reddi, Manzil Zaheer, Suvrit Sra, Barnabas Poczos, Francis Bach, Ruslan Salakhutdinov, and Alexander J Smola. A generic approach for escaping saddle points. *ArXiv e-prints*, abs/1709.01434, September 2017.

[41] Mark Schmidt, Nicolas Le Roux, and Francis Bach. Minimizing finite sums with the stochastic average gradient. *ArXiv e-prints*, abs/1309.2388, September 2013. Preliminary version appeared in NIPS 2012.

[42] Nicol N Schraudolph. Fast curvature matrix-vector products for second-order gradient descent. *Neural computation*, 14(7):1723–1738, 2002.

[43] Shai Shalev-Shwartz. Online Learning and Online Convex Optimization. *Foundations and Trends in Machine Learning*, 4(2):107–194, 2012. ISSN 1935-8237.

[44] Shai Shalev-Shwartz. SDCA without Duality, Regularization, and Individual Convexity. In *ICML*, 2016.

[45] Ruoyu Sun and Zhi-Quan Luo. Guaranteed Matrix Completion via Nonconvex Factorization. In *FOCS*, 2015.

[46] Nilesh Tripuraneni, Mitchell Stern, Chi Jin, Jeffrey Regier, and Michael I Jordan. Stochastic Cubic Regularization for Fast Nonconvex Optimization. *ArXiv e-prints*, abs/1711.02838, November 2017.

[47] Lin Xiao and Tong Zhang. A Proximal Stochastic Gradient Method with Progressive Variance Reduction. *SIAM Journal on Optimization*, 24(4):2057—-2075, 2014.

[48] Yi Xu and Tianbao Yang. First-order Stochastic Algorithms for Escaping From Saddle Points in Almost

Linear Time. *ArXiv e-prints*, abs/1711.01944, November 2017.

[49] Matthew D Zeiler. ADADELTA: an adaptive learning rate method. *ArXiv e-prints*, abs/1212.5701, 12 2012.

[50] Lijun Zhang, Mehrdad Mahdavi, and Rong Jin. Linear convergence with condition number independent access of full gradients. In *Advances in Neural Information Processing Systems*, pages 980–988, 2013.